\input amstex 
%
\catcode`\@=11

\def\East#1#2{\setboxz@h{$\m@th\ssize\;{#1}\;\;$}%
 \setbox@ne\hbox{$\m@th\ssize\;{#2}\;\;$}\setbox\tw@\hbox{$\m@th#2$}%
 \dimen@\minaw@
 \ifdim\wdz@>\dimen@ \dimen@\wdz@ \fi  \ifdim\wd@ne>\dimen@ \dimen@\wd@ne \fi
 \ifdim\wd\tw@>\z@
  \mathrel{\mathop{\hbox to\dimen@{\rightarrowfill}}\limits^{#1}_{#2}}%
 \else
  \mathrel{\mathop{\hbox to\dimen@{\rightarrowfill}}\limits^{#1}}%
 \fi}
\def\West#1#2{\setboxz@h{$\m@th\ssize\;\;{#1}\;$}%
 \setbox@ne\hbox{$\m@th\ssize\;\;{#2}\;$}\setbox\tw@\hbox{$\m@th#2$}%
 \dimen@\minaw@
 \ifdim\wdz@>\dimen@ \dimen@\wdz@ \fi \ifdim\wd@ne>\dimen@ \dimen@\wd@ne \fi
 \ifdim\wd\tw@>\z@
  \mathrel{\mathop{\hbox to\dimen@{\leftarrowfill}}\limits^{#1}_{#2}}%
 \else
  \mathrel{\mathop{\hbox to\dimen@{\leftarrowfill}}\limits^{#1}}%
 \fi}
\font\arrow@i=lams1
\font\arrow@ii=lams2
\font\arrow@iii=lams3
\font\arrow@iv=lams4
\font\arrow@v=lams5
\newbox\zer@
\newdimen\standardcgap
\standardcgap=40\p@
\newdimen\hunit
\hunit=\tw@\p@
\newdimen\standardrgap
\standardrgap=32\p@
\newdimen\vunit
\vunit=1.6\p@
\def\Cgaps#1{\RIfM@
  \standardcgap=#1\standardcgap\relax \hunit=#1\hunit\relax
 \else \nonmatherr@\Cgaps \fi}
\def\Rgaps#1{\RIfM@
  \standardrgap=#1\standardrgap\relax \vunit=#1\vunit\relax
 \else \nonmatherr@\Rgaps \fi}
\newdimen\getdim@
\def\getcgap@#1{\ifcase#1\or\getdim@\z@\else\getdim@\standardcgap\fi}
\def\getrgap@#1{\ifcase#1\getdim@\z@\else\getdim@\standardrgap\fi}
\def\cgaps#1{\RIfM@
 \cgaps@{#1}\edef\getcgap@##1{\i@=##1\relax\the\toks@}\toks@{}\else
 \nonmatherr@\cgaps\fi}
\def\rgaps#1{\RIfM@
 \rgaps@{#1}\edef\getrgap@##1{\i@=##1\relax\the\toks@}\toks@{}\else
 \nonmatherr@\rgaps\fi}
\def\Gaps@@{\gaps@@}
\def\cgaps@#1{\toks@{\ifcase\i@\or\getdim@=\z@}%
 \gaps@@\standardcgap#1;\gaps@@\gaps@@
 \edef\next@{\the\toks@\noexpand\else\noexpand\getdim@\noexpand\standardcgap
  \noexpand\fi}%
 \toks@=\expandafter{\next@}}
\def\rgaps@#1{\toks@{\ifcase\i@\getdim@=\z@}%
 \gaps@@\standardrgap#1;\gaps@@\gaps@@
 \edef\next@{\the\toks@\noexpand\else\noexpand\getdim@\noexpand\standardrgap
  \noexpand\fi}%
 \toks@=\expandafter{\next@}}
\def\gaps@@#1#2;#3{\mgaps@#1#2\mgaps@
 \edef\next@{\the\toks@\noexpand\or\noexpand\getdim@
  \noexpand#1\the\mgapstoks@@}%
 \global\toks@=\expandafter{\next@}%
 \DN@{#3}%
 \ifx\next@\Gaps@@\gdef\next@##1\gaps@@{}\else
  \gdef\next@{\gaps@@#1#3}\fi\next@}
\def\mgaps@#1{\let\mgapsnext@#1\FN@\mgaps@@}
\def\mgaps@@{\ifx\next\space@\DN@. {\FN@\mgaps@@}\else
 \DN@.{\FN@\mgaps@@@}\fi\next@.}
\def\mgaps@@@{\ifx\next\w\let\next@\mgaps@@@@\else
 \let\next@\mgaps@@@@@\fi\next@}
\newtoks\mgapstoks@@
\def\mgaps@@@@@#1\mgaps@{\getdim@\mgapsnext@\getdim@#1\getdim@
 \edef\next@{\noexpand\getdim@\the\getdim@}%
 \mgapstoks@@=\expandafter{\next@}}
\def\mgaps@@@@\w#1#2\mgaps@{\mgaps@@@@@#2\mgaps@
 \setbox\zer@\hbox{$\m@th\hskip15\p@\tsize@#1$}%
 \dimen@\wd\zer@
 \ifdim\dimen@>\getdim@ \getdim@\dimen@ \fi
 \edef\next@{\noexpand\getdim@\the\getdim@}%
 \mgapstoks@@=\expandafter{\next@}}
\def\changewidth#1#2{\setbox\zer@\hbox{$\m@th#2}%
 \hbox to\wd\zer@{\hss$\m@th#1$\hss}}
\atdef@({\FN@\ARROW@}
\def\ARROW@{\ifx\next)\let\next@\OPTIONS@\else
 \DN@{\csname\string @(\endcsname}\fi\next@}
\newif\ifoptions@
\def\OPTIONS@){\ifoptions@\let\next@\relax\else
 \DN@{\options@true\begingroup\optioncodes@}\fi\next@}
\newif\ifN@
\newif\ifE@
\newif\ifNESW@
\newif\ifH@
\newif\ifV@
\newif\ifHshort@
\expandafter\def\csname\string @(\endcsname #1,#2){%
 \ifoptions@\let\next@\endgroup\else\let\next@\relax\fi\next@
 \N@false\E@false\H@false\V@false\Hshort@false
 \ifnum#1>\z@\E@true\fi
 \ifnum#1=\z@\V@true\tX@false\tY@false\a@false\fi
 \ifnum#2>\z@\N@true\fi
 \ifnum#2=\z@\H@true\tX@false\tY@false\a@false\ifshort@\Hshort@true\fi\fi
 \NESW@false
 \ifN@\ifE@\NESW@true\fi\else\ifE@\else\NESW@true\fi\fi
 \arrow@{#1}{#2}%
 \global\options@false
 \global\scount@\z@\global\tcount@\z@\global\arrcount@\z@
 \global\s@false\global\sxdimen@\z@\global\sydimen@\z@
 \global\tX@false\global\tXdimen@i\z@\global\tXdimen@ii\z@
 \global\tY@false\global\tYdimen@i\z@\global\tYdimen@ii\z@
 \global\a@false\global\exacount@\z@
 \global\x@false\global\xdimen@\z@
 \global\X@false\global\Xdimen@\z@
 \global\y@false\global\ydimen@\z@
 \global\Y@false\global\Ydimen@\z@
 \global\p@false\global\pdimen@\z@
 \global\label@ifalse\global\label@iifalse
 \global\dl@ifalse\global\ldimen@i\z@
 \global\dl@iifalse\global\ldimen@ii\z@
 \global\short@false\global\unshort@false}
\newif\iflabel@i
\newif\iflabel@ii
\newcount\scount@
\newcount\tcount@
\newcount\arrcount@
\newif\ifs@
\newdimen\sxdimen@
\newdimen\sydimen@
\newif\iftX@
\newdimen\tXdimen@i
\newdimen\tXdimen@ii
\newif\iftY@
\newdimen\tYdimen@i
\newdimen\tYdimen@ii
\newif\ifa@
\newcount\exacount@
\newif\ifx@
\newdimen\xdimen@
\newif\ifX@
\newdimen\Xdimen@
\newif\ify@
\newdimen\ydimen@
\newif\ifY@
\newdimen\Ydimen@
\newif\ifp@
\newdimen\pdimen@
\newif\ifdl@i
\newif\ifdl@ii
\newdimen\ldimen@i
\newdimen\ldimen@ii
\newif\ifshort@
\newif\ifunshort@
\def\zero@#1{\ifnum\scount@=\z@
 \if#1e\global\scount@\m@ne\else
 \if#1t\global\scount@\tw@\else
 \if#1h\global\scount@\thr@@\else
 \if#1'\global\scount@6 \else
 \if#1`\global\scount@7 \else
 \if#1(\global\scount@8 \else
 \if#1)\global\scount@9 \else
 \if#1s\global\scount@12 \else
 \if#1H\global\scount@13 \else
 \Err@{\Invalid@@ option \string\0}\fi\fi\fi\fi\fi\fi\fi\fi\fi
 \fi}
\def\one@#1{\ifnum\tcount@=\z@
 \if#1e\global\tcount@\m@ne\else
 \if#1h\global\tcount@\tw@\else
 \if#1t\global\tcount@\thr@@\else
 \if#1'\global\tcount@4 \else
 \if#1`\global\tcount@5 \else
 \if#1(\global\tcount@10 \else
 \if#1)\global\tcount@11 \else
 \if#1s\global\tcount@12 \else
 \if#1H\global\tcount@13 \else
 \Err@{\Invalid@@ option \string\1}\fi\fi\fi\fi\fi\fi\fi\fi\fi
 \fi}
\def\a@#1{\ifnum\arrcount@=\z@
 \if#10\global\arrcount@\m@ne\else
 \if#1+\global\arrcount@\@ne\else
 \if#1-\global\arrcount@\tw@\else
 \if#1=\global\arrcount@\thr@@\else
 \Err@{\Invalid@@ option \string\a}\fi\fi\fi\fi
 \fi}
\def\ds@(#1;#2){\ifs@\else
 \global\s@true
 \sxdimen@\hunit \global\sxdimen@#1\sxdimen@\relax
 \sydimen@\vunit \global\sydimen@#2\sydimen@\relax
 \fi}
\def\dtX@(#1;#2){\iftX@\else
 \global\tX@true
 \tXdimen@i\hunit \global\tXdimen@i#1\tXdimen@i\relax
 \tXdimen@ii\vunit \global\tXdimen@ii#2\tXdimen@ii\relax
 \fi}
\def\dtY@(#1;#2){\iftY@\else
 \global\tY@true
 \tYdimen@i\hunit \global\tYdimen@i#1\tYdimen@i\relax
 \tYdimen@ii\vunit \global\tYdimen@ii#2\tYdimen@ii\relax
 \fi}
\def\da@#1{\ifa@\else\global\a@true\global\exacount@#1\relax\fi}
\def\dx@#1{\ifx@\else
 \global\x@true
 \xdimen@\hunit \global\xdimen@#1\xdimen@\relax
 \fi}
\def\dX@#1{\ifX@\else
 \global\X@true
 \Xdimen@\hunit \global\Xdimen@#1\Xdimen@\relax
 \fi}
\def\dy@#1{\ify@\else
 \global\y@true
 \ydimen@\vunit \global\ydimen@#1\ydimen@\relax
 \fi}
\def\dY@#1{\ifY@\else
 \global\Y@true
 \Ydimen@\vunit \global\Ydimen@#1\Ydimen@\relax
 \fi}
\def\p@@#1{\ifp@\else
 \global\p@true
 \pdimen@\hunit \divide\pdimen@\tw@ \global\pdimen@#1\pdimen@\relax
 \fi}
\def\L@#1{\iflabel@i\else
 \global\label@itrue \gdef\label@i{#1}%
 \fi}
\def\l@#1{\iflabel@ii\else
 \global\label@iitrue \gdef\label@ii{#1}%
 \fi}
\def\dL@#1{\ifdl@i\else
 \global\dl@itrue \ldimen@i\hunit \global\ldimen@i#1\ldimen@i\relax
 \fi}
\def\dl@#1{\ifdl@ii\else
 \global\dl@iitrue \ldimen@ii\hunit \global\ldimen@ii#1\ldimen@ii\relax
 \fi}
\def\s@{\ifunshort@\else\global\short@true\fi}
\def\uns@{\ifshort@\else\global\unshort@true\global\short@false\fi}
\def\optioncodes@{\let\0\zero@\let\1\one@\let\a\a@\let\ds\ds@\let\dtX\dtX@
 \let\dtY\dtY@\let\da\da@\let\dx\dx@\let\dX\dX@\let\dY\dY@\let\dy\dy@
 \let\p\p@@\let\L\L@\let\l\l@\let\dL\dL@\let\dl\dl@\let\s\s@\let\uns\uns@}
\def\slopes@{\\161\\152\\143\\134\\255\\126\\357\\238\\349\\45{10}\\56{11}%
 \\11{12}\\65{13}\\54{14}\\43{15}\\32{16}\\53{17}\\21{18}\\52{19}\\31{20}%
 \\41{21}\\51{22}\\61{23}}
\newcount\tan@i
\newcount\tan@ip
\newcount\tan@ii
\newcount\tan@iip
\newdimen\slope@i
\newdimen\slope@ip
\newdimen\slope@ii
\newdimen\slope@iip
\newcount\angcount@
\newcount\extracount@
\def\slope@{{\slope@i=\secondy@ \advance\slope@i-\firsty@
 \ifN@\else\multiply\slope@i\m@ne\fi
 \slope@ii=\secondx@ \advance\slope@ii-\firstx@
 \ifE@\else\multiply\slope@ii\m@ne\fi
 \ifdim\slope@ii<\z@
  \global\tan@i6 \global\tan@ii\@ne \global\angcount@23
 \else
  \dimen@\slope@i \multiply\dimen@6
  \ifdim\dimen@<\slope@ii
   \global\tan@i\@ne \global\tan@ii6 \global\angcount@\@ne
  \else
   \dimen@\slope@ii \multiply\dimen@6
   \ifdim\dimen@<\slope@i
    \global\tan@i6 \global\tan@ii\@ne \global\angcount@23
   \else
    \tan@ip\z@ \tan@iip \@ne
    \def\\##1##2##3{\global\angcount@=##3\relax
     \slope@ip\slope@i \slope@iip\slope@ii
     \multiply\slope@iip##1\relax \multiply\slope@ip##2\relax
     \ifdim\slope@iip<\slope@ip
      \global\tan@ip=##1\relax \global\tan@iip=##2\relax
     \else
      \global\tan@i=##1\relax \global\tan@ii=##2\relax
      \def\\####1####2####3{}%
     \fi}%
    \slopes@
    \slope@i=\secondy@ \advance\slope@i-\firsty@
    \ifN@\else\multiply\slope@i\m@ne\fi
    \multiply\slope@i\tan@ii \multiply\slope@i\tan@iip \multiply\slope@i\tw@
    \count@\tan@i \multiply\count@\tan@iip
    \extracount@\tan@ip \multiply\extracount@\tan@ii
    \advance\count@\extracount@
    \slope@ii=\secondx@ \advance\slope@ii-\firstx@
    \ifE@\else\multiply\slope@ii\m@ne\fi
    \multiply\slope@ii\count@
    \ifdim\slope@i<\slope@ii
     \global\tan@i=\tan@ip \global\tan@ii=\tan@iip
     \global\advance\angcount@\m@ne
    \fi
   \fi
  \fi
 \fi}%
}
\def\slope@a#1{{\def\\##1##2##3{\ifnum##3=#1\global\tan@i=##1\relax
 \global\tan@ii=##2\relax\fi}\slopes@}}
\newcount\i@
\newcount\j@
\newcount\colcount@
\newcount\Colcount@
\newcount\tcolcount@
\newdimen\rowht@
\newdimen\rowdp@
\newcount\rowcount@
\newcount\Rowcount@
\newcount\maxcolrow@
\newtoks\colwidthtoks@
\newtoks\Rowheighttoks@
\newtoks\Rowdepthtoks@
\newtoks\widthtoks@
\newtoks\Widthtoks@
\newtoks\heighttoks@
\newtoks\Heighttoks@
\newtoks\depthtoks@
\newtoks\Depthtoks@
\newif\iffirstnewCDcr@
\def\dotoks@i{%
 \global\widthtoks@=\expandafter{\the\widthtoks@\else\getdim@\z@\fi}%
 \global\heighttoks@=\expandafter{\the\heighttoks@\else\getdim@\z@\fi}%
 \global\depthtoks@=\expandafter{\the\depthtoks@\else\getdim@\z@\fi}}
\def\dotoks@ii{%
 \global\widthtoks@{\ifcase\j@}%
 \global\heighttoks@{\ifcase\j@}%
 \global\depthtoks@{\ifcase\j@}}
\def\prenewCD@#1\endnewCD{\setbox\zer@
 \vbox{%
  \def\arrow@##1##2{{}}%
  \rowcount@\m@ne \colcount@\z@ \Colcount@\z@
  \firstnewCDcr@true \toks@{}%
  \widthtoks@{\ifcase\j@}%
  \Widthtoks@{\ifcase\i@}%
  \heighttoks@{\ifcase\j@}%
  \Heighttoks@{\ifcase\i@}%
  \depthtoks@{\ifcase\j@}%
  \Depthtoks@{\ifcase\i@}%
  \Rowheighttoks@{\ifcase\i@}%
  \Rowdepthtoks@{\ifcase\i@}%
  \Let@
  \everycr{%
   \noalign{%
    \global\advance\rowcount@\@ne
    \ifnum\colcount@<\Colcount@
    \else
     \global\Colcount@=\colcount@ \global\maxcolrow@=\rowcount@
    \fi
    \global\colcount@\z@
    \iffirstnewCDcr@
     \global\firstnewCDcr@false
    \else
     \edef\next@{\the\Rowheighttoks@\noexpand\or\noexpand\getdim@\the\rowht@}%
      \global\Rowheighttoks@=\expandafter{\next@}%
     \edef\next@{\the\Rowdepthtoks@\noexpand\or\noexpand\getdim@\the\rowdp@}%
      \global\Rowdepthtoks@=\expandafter{\next@}%
     \global\rowht@\z@ \global\rowdp@\z@
     \dotoks@i
     \edef\next@{\the\Widthtoks@\noexpand\or\the\widthtoks@}%
      \global\Widthtoks@=\expandafter{\next@}%
     \edef\next@{\the\Heighttoks@\noexpand\or\the\heighttoks@}%
      \global\Heighttoks@=\expandafter{\next@}%
     \edef\next@{\the\Depthtoks@\noexpand\or\the\depthtoks@}%
      \global\Depthtoks@=\expandafter{\next@}%
     \dotoks@ii
    \fi}}%
  \tabskip\z@
  \halign{&\setbox\zer@\hbox{\vrule height10\p@ width\z@ depth\z@
   $\m@th\displaystyle{##}$}\copy\zer@
   \ifdim\ht\zer@>\rowht@ \global\rowht@\ht\zer@ \fi
   \ifdim\dp\zer@>\rowdp@ \global\rowdp@\dp\zer@ \fi
   \global\advance\colcount@\@ne
   \edef\next@{\the\widthtoks@\noexpand\or\noexpand\getdim@\the\wd\zer@}%
    \global\widthtoks@=\expandafter{\next@}%
   \edef\next@{\the\heighttoks@\noexpand\or\noexpand\getdim@\the\ht\zer@}%
    \global\heighttoks@=\expandafter{\next@}%
   \edef\next@{\the\depthtoks@\noexpand\or\noexpand\getdim@\the\dp\zer@}%
    \global\depthtoks@=\expandafter{\next@}%
   \cr#1\crcr}}%
 \Rowcount@=\rowcount@
 \global\Widthtoks@=\expandafter{\the\Widthtoks@\fi\relax}%
 \edef\Width@##1##2{\i@=##1\relax\j@=##2\relax\the\Widthtoks@}%
 \global\Heighttoks@=\expandafter{\the\Heighttoks@\fi\relax}%
 \edef\Height@##1##2{\i@=##1\relax\j@=##2\relax\the\Heighttoks@}%
 \global\Depthtoks@=\expandafter{\the\Depthtoks@\fi\relax}%
 \edef\Depth@##1##2{\i@=##1\relax\j@=##2\relax\the\Depthtoks@}%
 \edef\next@{\the\Rowheighttoks@\noexpand\fi\relax}%
 \global\Rowheighttoks@=\expandafter{\next@}%
 \edef\Rowheight@##1{\i@=##1\relax\the\Rowheighttoks@}%
 \edef\next@{\the\Rowdepthtoks@\noexpand\fi\relax}%
 \global\Rowdepthtoks@=\expandafter{\next@}%
 \edef\Rowdepth@##1{\i@=##1\relax\the\Rowdepthtoks@}%
 \colwidthtoks@{\fi}%
 \setbox\zer@\vbox{%
  \unvbox\zer@
  \count@\rowcount@
  \loop
   \unskip\unpenalty
   \setbox\zer@\lastbox
   \ifnum\count@>\maxcolrow@ \advance\count@\m@ne
   \repeat
  \hbox{%
   \unhbox\zer@
   \count@\z@
   \loop
    \unskip
    \setbox\zer@\lastbox
    \edef\next@{\noexpand\or\noexpand\getdim@\the\wd\zer@\the\colwidthtoks@}%
     \global\colwidthtoks@=\expandafter{\next@}%
    \advance\count@\@ne
    \ifnum\count@<\Colcount@
    \repeat}}%
 \edef\next@{\noexpand\ifcase\noexpand\i@\the\colwidthtoks@}%
  \global\colwidthtoks@=\expandafter{\next@}%
 \edef\Colwidth@##1{\i@=##1\relax\the\colwidthtoks@}%
 \colwidthtoks@{}\Rowheighttoks@{}\Rowdepthtoks@{}\widthtoks@{}%
 \Widthtoks@{}\heighttoks@{}\Heighttoks@{}\depthtoks@{}\Depthtoks@{}%
}
\newcount\xoff@
\newcount\yoff@
\newcount\endcount@
\newcount\rcount@
\newdimen\firstx@
\newdimen\firsty@
\newdimen\secondx@
\newdimen\secondy@
\newdimen\tocenter@
\newdimen\charht@
\newdimen\charwd@
\def\outside@{\Err@{This arrow points outside the \string\newCD}}
\newif\ifsvertex@
\newif\iftvertex@
\def\arrow@#1#2{\xoff@=#1\relax\yoff@=#2\relax
 \count@\rowcount@ \advance\count@-\yoff@
 \ifnum\count@<\@ne \outside@ \else \ifnum\count@>\Rowcount@ \outside@ \fi\fi
 \count@\colcount@ \advance\count@\xoff@
 \ifnum\count@<\@ne \outside@ \else \ifnum\count@>\Colcount@ \outside@\fi\fi
 \tcolcount@\colcount@ \advance\tcolcount@\xoff@
 \Width@\rowcount@\colcount@ \tocenter@=-\getdim@ \divide\tocenter@\tw@
 \ifdim\getdim@=\z@
  \firstx@\z@ \firsty@\mathaxis@ \svertex@true
 \else
  \svertex@false
  \ifHshort@
   \Colwidth@\colcount@
    \ifE@ \firstx@=.5\getdim@ \else \firstx@=-.5\getdim@ \fi
  \else
   \ifE@ \firstx@=\getdim@ \else \firstx@=-\getdim@ \fi
   \divide\firstx@\tw@
  \fi
  \ifE@
   \ifH@ \advance\firstx@\thr@@\p@ \else \advance\firstx@-\thr@@\p@ \fi
  \else
   \ifH@ \advance\firstx@-\thr@@\p@ \else \advance\firstx@\thr@@\p@ \fi
  \fi
  \ifN@
   \Height@\rowcount@\colcount@ \firsty@=\getdim@
   \ifV@ \advance\firsty@\thr@@\p@ \fi
  \else
   \ifV@
    \Depth@\rowcount@\colcount@ \firsty@=-\getdim@
    \advance\firsty@-\thr@@\p@
   \else
    \firsty@\z@
   \fi
  \fi
 \fi
 \ifV@
 \else
  \Colwidth@\colcount@
  \ifE@ \secondx@=\getdim@ \else \secondx@=-\getdim@ \fi
  \divide\secondx@\tw@
  \ifE@ \else \getcgap@\colcount@ \advance\secondx@-\getdim@ \fi
  \endcount@=\colcount@ \advance\endcount@\xoff@
  \count@=\colcount@
  \ifE@
   \advance\count@\@ne
   \loop
    \ifnum\count@<\endcount@
    \Colwidth@\count@ \advance\secondx@\getdim@
    \getcgap@\count@ \advance\secondx@\getdim@
    \advance\count@\@ne
    \repeat
  \else
   \advance\count@\m@ne
   \loop
    \ifnum\count@>\endcount@
    \Colwidth@\count@ \advance\secondx@-\getdim@
    \getcgap@\count@ \advance\secondx@-\getdim@
    \advance\count@\m@ne
    \repeat
  \fi
  \Colwidth@\count@ \divide\getdim@\tw@
  \ifHshort@
  \else
   \ifE@ \advance\secondx@\getdim@ \else \advance\secondx@-\getdim@ \fi
  \fi
  \ifE@ \getcgap@\count@ \advance\secondx@\getdim@ \fi
  \rcount@\rowcount@ \advance\rcount@-\yoff@
  \Width@\rcount@\count@ \divide\getdim@\tw@
  \tvertex@false
  \ifH@\ifdim\getdim@=\z@\tvertex@true\Hshort@false\fi\fi
  \ifHshort@
  \else
   \ifE@ \advance\secondx@-\getdim@ \else \advance\secondx@\getdim@ \fi
  \fi
  \iftvertex@
   \advance\secondx@.4\p@
  \else
   \ifE@ \advance\secondx@-\thr@@\p@ \else \advance\secondx@\thr@@\p@ \fi
  \fi
 \fi
 \ifH@
 \else
  \ifN@
   \Rowheight@\rowcount@ \secondy@\getdim@
  \else
   \Rowdepth@\rowcount@ \secondy@-\getdim@
   \getrgap@\rowcount@ \advance\secondy@-\getdim@
  \fi
  \endcount@=\rowcount@ \advance\endcount@-\yoff@
  \count@=\rowcount@
  \ifN@
   \advance\count@\m@ne
   \loop
    \ifnum\count@>\endcount@
    \Rowheight@\count@ \advance\secondy@\getdim@
    \Rowdepth@\count@ \advance\secondy@\getdim@
    \getrgap@\count@ \advance\secondy@\getdim@
    \advance\count@\m@ne
    \repeat
  \else
   \advance\count@\@ne
   \loop
    \ifnum\count@<\endcount@
    \Rowheight@\count@ \advance\secondy@-\getdim@
    \Rowdepth@\count@ \advance\secondy@-\getdim@
    \getrgap@\count@ \advance\secondy@-\getdim@
    \advance\count@\@ne
    \repeat
  \fi
  \tvertex@false
  \ifV@\Width@\count@\colcount@\ifdim\getdim@=\z@\tvertex@true\fi\fi
  \ifN@
   \getrgap@\count@ \advance\secondy@\getdim@
   \Rowdepth@\count@ \advance\secondy@\getdim@
   \iftvertex@
    \advance\secondy@\mathaxis@
   \else
    \Depth@\count@\tcolcount@ \advance\secondy@-\getdim@
    \advance\secondy@-\thr@@\p@
   \fi
  \else
   \Rowheight@\count@ \advance\secondy@-\getdim@
   \iftvertex@
    \advance\secondy@\mathaxis@
   \else
    \Height@\count@\tcolcount@ \advance\secondy@\getdim@
    \advance\secondy@\thr@@\p@
   \fi
  \fi
 \fi
 \ifV@\else\advance\firstx@\sxdimen@\fi
 \ifH@\else\advance\firsty@\sydimen@\fi
 \iftX@
  \advance\secondy@\tXdimen@ii
  \advance\secondx@\tXdimen@i
  \slope@
 \else
  \iftY@
   \advance\secondy@\tYdimen@ii
   \advance\secondx@\tYdimen@i
   \slope@
   \secondy@=\secondx@ \advance\secondy@-\firstx@
   \ifNESW@ \else \multiply\secondy@\m@ne \fi
   \multiply\secondy@\tan@i \divide\secondy@\tan@ii \advance\secondy@\firsty@
  \else
   \ifa@
    \slope@
    \ifNESW@ \global\advance\angcount@\exacount@ \else
      \global\advance\angcount@-\exacount@ \fi
    \ifnum\angcount@>23 \angcount@23 \fi
    \ifnum\angcount@<\@ne \angcount@\@ne \fi
    \slope@a\angcount@
    \ifY@
     \advance\secondy@\Ydimen@
    \else
     \ifX@
      \advance\secondx@\Xdimen@
      \dimen@\secondx@ \advance\dimen@-\firstx@
      \ifNESW@\else\multiply\dimen@\m@ne\fi
      \multiply\dimen@\tan@i \divide\dimen@\tan@ii
      \advance\dimen@\firsty@ \secondy@=\dimen@
     \fi
    \fi
   \else
    \ifH@\else\ifV@\else\slope@\fi\fi
   \fi
  \fi
 \fi
 \ifH@\else\ifV@\else\ifsvertex@\else
  \dimen@=6\p@ \multiply\dimen@\tan@ii
  \count@=\tan@i \advance\count@\tan@ii \divide\dimen@\count@
  \ifE@ \advance\firstx@\dimen@ \else \advance\firstx@-\dimen@ \fi
  \multiply\dimen@\tan@i \divide\dimen@\tan@ii
  \ifN@ \advance\firsty@\dimen@ \else \advance\firsty@-\dimen@ \fi
 \fi\fi\fi
 \ifp@
  \ifH@\else\ifV@\else
   \getcos@\pdimen@ \advance\firsty@\dimen@ \advance\secondy@\dimen@
   \ifNESW@ \advance\firstx@-\dimen@ii \else \advance\firstx@\dimen@ii \fi
  \fi\fi
 \fi
 \ifH@\else\ifV@\else
  \ifnum\tan@i>\tan@ii
   \charht@=10\p@ \charwd@=10\p@
   \multiply\charwd@\tan@ii \divide\charwd@\tan@i
  \else
   \charwd@=10\p@ \charht@=10\p@
   \divide\charht@\tan@ii \multiply\charht@\tan@i
  \fi
  \ifnum\tcount@=\thr@@
   \ifN@ \advance\secondy@-.3\charht@ \else\advance\secondy@.3\charht@ \fi
  \fi
  \ifnum\scount@=\tw@
   \ifE@ \advance\firstx@.3\charht@ \else \advance\firstx@-.3\charht@ \fi
  \fi
  \ifnum\tcount@=12
   \ifN@ \advance\secondy@-\charht@ \else \advance\secondy@\charht@ \fi
  \fi
  \iftY@
  \else
   \ifa@
    \ifX@
    \else
     \secondx@\secondy@ \advance\secondx@-\firsty@
     \ifNESW@\else\multiply\secondx@\m@ne\fi
     \multiply\secondx@\tan@ii \divide\secondx@\tan@i
     \advance\secondx@\firstx@
    \fi
   \fi
  \fi
 \fi\fi
 \ifH@\harrow@\else\ifV@\varrow@\else\arrow@@\fi\fi}
\newdimen\mathaxis@
\mathaxis@90\p@ \divide\mathaxis@36
\def\harrow@b{\ifE@\hskip\tocenter@\hskip\firstx@\fi}
\def\harrow@bb{\ifE@\hskip\xdimen@\else\hskip\Xdimen@\fi}
\def\harrow@e{\ifE@\else\hskip-\firstx@\hskip-\tocenter@\fi}
\def\harrow@ee{\ifE@\hskip-\Xdimen@\else\hskip-\xdimen@\fi}
\def\harrow@{\dimen@\secondx@\advance\dimen@-\firstx@
 \ifE@ \let\next@\rlap \else  \multiply\dimen@\m@ne \let\next@\llap \fi
 \next@{%
  \harrow@b
  \smash{\raise\pdimen@\hbox to\dimen@
   {\harrow@bb\arrow@ii
    \ifnum\arrcount@=\m@ne \else \ifnum\arrcount@=\thr@@ \else
     \ifE@
      \ifnum\scount@=\m@ne
      \else
       \ifcase\scount@\or\or\char118 \or\char117 \or\or\or\char119 \or
       \char120 \or\char121 \or\char122 \or\or\or\arrow@i\char125 \or
       \char117 \hskip\thr@@\p@\char117 \hskip-\thr@@\p@\fi
      \fi
     \else
      \ifnum\tcount@=\m@ne
      \else
       \ifcase\tcount@\char117 \or\or\char117 \or\char118 \or\char119 \or
       \char120\or\or\or\or\or\char121 \or\char122 \or\arrow@i\char125
       \or\char117 \hskip\thr@@\p@\char117 \hskip-\thr@@\p@\fi
      \fi
     \fi
    \fi\fi
    \dimen@\mathaxis@ \advance\dimen@.2\p@
    \dimen@ii\mathaxis@ \advance\dimen@ii-.2\p@
    \ifnum\arrcount@=\m@ne
     \let\leads@\null
    \else
     \ifcase\arrcount@
      \def\leads@{\hrule height\dimen@ depth-\dimen@ii}\or
      \def\leads@{\hrule height\dimen@ depth-\dimen@ii}\or
      \def\leads@{\hbox to10\p@{%
       \leaders\hrule height\dimen@ depth-\dimen@ii\hfil
       \hfil
      \leaders\hrule height\dimen@ depth-\dimen@ii\hskip\z@ plus2fil\relax
       \hfil
       \leaders\hrule height\dimen@ depth-\dimen@ii\hfil}}\or
     \def\leads@{\hbox{\hbox to10\p@{\dimen@\mathaxis@ \advance\dimen@1.2\p@
       \dimen@ii\dimen@ \advance\dimen@ii-.4\p@
       \leaders\hrule height\dimen@ depth-\dimen@ii\hfil}%
       \kern-10\p@
       \hbox to10\p@{\dimen@\mathaxis@ \advance\dimen@-1.2\p@
       \dimen@ii\dimen@ \advance\dimen@ii-.4\p@
       \leaders\hrule height\dimen@ depth-\dimen@ii\hfil}}}\fi
    \fi
    \cleaders\leads@\hfil
    \ifnum\arrcount@=\m@ne\else\ifnum\arrcount@=\thr@@\else
     \arrow@i
     \ifE@
      \ifnum\tcount@=\m@ne
      \else
       \ifcase\tcount@\char119 \or\or\char119 \or\char120 \or\char121 \or
       \char122 \or \or\or\or\or\char123\or\char124 \or
       \char125 \or\char119 \hskip-\thr@@\p@\char119 \hskip\thr@@\p@\fi
      \fi
     \else
      \ifcase\scount@\or\or\char120 \or\char119 \or\or\or\char121 \or\char122
      \or\char123 \or\char124 \or\or\or\char125 \or
      \char119 \hskip-\thr@@\p@\char119 \hskip\thr@@\p@\fi
     \fi
    \fi\fi
    \harrow@ee}}%
  \harrow@e}%
 \iflabel@i
  \dimen@ii\z@ \setbox\zer@\hbox{$\m@th\tsize@@\label@i$}%
  \ifnum\arrcount@=\m@ne
  \else
   \advance\dimen@ii\mathaxis@
   \advance\dimen@ii\dp\zer@ \advance\dimen@ii\tw@\p@
   \ifnum\arrcount@=\thr@@ \advance\dimen@ii\tw@\p@ \fi
  \fi
  \advance\dimen@ii\pdimen@
  \next@{\harrow@b\smash{\raise\dimen@ii\hbox to\dimen@
   {\harrow@bb\hskip\tw@\ldimen@i\hfil\box\zer@\hfil\harrow@ee}}\harrow@e}%
 \fi
 \iflabel@ii
  \ifnum\arrcount@=\m@ne
  \else
   \setbox\zer@\hbox{$\m@th\tsize@\label@ii$}%
   \dimen@ii-\ht\zer@ \advance\dimen@ii-\tw@\p@
   \ifnum\arrcount@=\thr@@ \advance\dimen@ii-\tw@\p@ \fi
   \advance\dimen@ii\mathaxis@ \advance\dimen@ii\pdimen@
   \next@{\harrow@b\smash{\raise\dimen@ii\hbox to\dimen@
    {\harrow@bb\hskip\tw@\ldimen@ii\hfil\box\zer@\hfil\harrow@ee}}\harrow@e}%
  \fi
 \fi}
\let\tsize@\tsize
\def\tsizenewCDlabels{\let\tsize@\tsize}
\def\ssizenewCDlabels{\let\tsize@\ssize}
\def\tsize@@{\ifnum\arrcount@=\m@ne\else\tsize@\fi}
\def\varrow@{\dimen@\secondy@ \advance\dimen@-\firsty@
 \ifN@ \else \multiply\dimen@\m@ne \fi
 \setbox\zer@\vbox to\dimen@
  {\ifN@ \vskip-\Ydimen@ \else \vskip\ydimen@ \fi
   \ifnum\arrcount@=\m@ne\else\ifnum\arrcount@=\thr@@\else
    \hbox{\arrow@iii
     \ifN@
      \ifnum\tcount@=\m@ne
      \else
       \ifcase\tcount@\char117 \or\or\char117 \or\char118 \or\char119 \or
       \char120 \or\or\or\or\or\char121 \or\char122 \or\char123 \or
       \vbox{\hbox{\char117 }\nointerlineskip\vskip\thr@@\p@
       \hbox{\char117 }\vskip-\thr@@\p@}\fi
      \fi
     \else
      \ifcase\scount@\or\or\char118 \or\char117 \or\or\or\char119 \or
      \char120 \or\char121 \or\char122 \or\or\or\char123 \or
      \vbox{\hbox{\char117 }\nointerlineskip\vskip\thr@@\p@
      \hbox{\char117 }\vskip-\thr@@\p@}\fi
     \fi}%
    \nointerlineskip
   \fi\fi
   \ifnum\arrcount@=\m@ne
    \let\leads@\null
   \else
    \ifcase\arrcount@\let\leads@\vrule\or\let\leads@\vrule\or
    \def\leads@{\vbox to10\p@{%
     \hrule height 1.67\p@ depth\z@ width.4\p@
     \vfil
     \hrule height 3.33\p@ depth\z@ width.4\p@
     \vfil
     \hrule height 1.67\p@ depth\z@ width.4\p@}}\or
    \def\leads@{\hbox{\vrule height\p@\hskip\tw@\p@\vrule}}\fi
   \fi
  \cleaders\leads@\vfill\nointerlineskip
   \ifnum\arrcount@=\m@ne\else\ifnum\arrcount@=\thr@@\else
    \hbox{\arrow@iv
     \ifN@
      \ifcase\scount@\or\or\char118 \or\char117 \or\or\or\char119 \or
      \char120 \or\char121 \or\char122 \or\or\or\arrow@iii\char123 \or
      \vbox{\hbox{\char117 }\nointerlineskip\vskip-\thr@@\p@
      \hbox{\char117 }\vskip\thr@@\p@}\fi
     \else
      \ifnum\tcount@=\m@ne
      \else
       \ifcase\tcount@\char117 \or\or\char117 \or\char118 \or\char119 \or
       \char120 \or\or\or\or\or\char121 \or\char122 \or\arrow@iii\char123 \or
       \vbox{\hbox{\char117 }\nointerlineskip\vskip-\thr@@\p@
       \hbox{\char117 }\vskip\thr@@\p@}\fi
      \fi
     \fi}%
   \fi\fi
   \ifN@\vskip\ydimen@\else\vskip-\Ydimen@\fi}%
 \ifN@
  \dimen@ii\firsty@
 \else
  \dimen@ii-\firsty@ \advance\dimen@ii\ht\zer@ \multiply\dimen@ii\m@ne
 \fi
 \rlap{\smash{\hskip\tocenter@ \hskip\pdimen@ \raise\dimen@ii \box\zer@}}%
 \iflabel@i
  \setbox\zer@\vbox to\dimen@{\vfil
   \hbox{$\m@th\tsize@@\label@i$}\vskip\tw@\ldimen@i\vfil}%
  \rlap{\smash{\hskip\tocenter@ \hskip\pdimen@
  \ifnum\arrcount@=\m@ne \let\next@\relax \else \let\next@\llap \fi
  \next@{\raise\dimen@ii\hbox{\ifnum\arrcount@=\m@ne \hskip-.5\wd\zer@ \fi
   \box\zer@ \ifnum\arrcount@=\m@ne \else \hskip\tw@\p@ \fi}}}}%
 \fi
 \iflabel@ii
  \ifnum\arrcount@=\m@ne
  \else
   \setbox\zer@\vbox to\dimen@{\vfil
    \hbox{$\m@th\tsize@\label@ii$}\vskip\tw@\ldimen@ii\vfil}%
   \rlap{\smash{\hskip\tocenter@ \hskip\pdimen@
   \rlap{\raise\dimen@ii\hbox{\ifnum\arrcount@=\thr@@ \hskip4.5\p@ \else
    \hskip2.5\p@ \fi\box\zer@}}}}%
  \fi
 \fi
}
\newdimen\goal@
\newdimen\shifted@
\newcount\Tcount@
\newcount\Scount@
\newbox\shaft@
\newcount\slcount@
\def\getcos@#1{%
 \ifnum\tan@i<\tan@ii
  \dimen@#1%
  \ifnum\slcount@<8 \count@9 \else \ifnum\slcount@<12 \count@8 \else
   \count@7 \fi\fi
  \multiply\dimen@\count@ \divide\dimen@10
  \dimen@ii\dimen@ \multiply\dimen@ii\tan@i \divide\dimen@ii\tan@ii
 \else
  \dimen@ii#1%
  \count@-\slcount@ \advance\count@24
  \ifnum\count@<8 \count@9 \else \ifnum\count@<12 \count@8
   \else\count@7 \fi\fi
  \multiply\dimen@ii\count@ \divide\dimen@ii10
  \dimen@\dimen@ii \multiply\dimen@\tan@ii \divide\dimen@\tan@i
 \fi}
\newdimen\adjust@
\def\Nnext@{\ifN@\let\next@\raise\else\let\next@\lower\fi}
\def\arrow@@{\slcount@\angcount@
 \ifNESW@
  \ifnum\angcount@<10
   \let\arrowfont@=\arrow@i \advance\angcount@\m@ne \multiply\angcount@13
  \else
   \ifnum\angcount@<19
    \let\arrowfont@=\arrow@ii \advance\angcount@-10 \multiply\angcount@13
   \else
    \let\arrowfont@=\arrow@iii \advance\angcount@-19 \multiply\angcount@13
  \fi\fi
  \Tcount@\angcount@
 \else
  \ifnum\angcount@<5
   \let\arrowfont@=\arrow@iii \advance\angcount@\m@ne \multiply\angcount@13
   \advance\angcount@65
  \else
   \ifnum\angcount@<14
    \let\arrowfont@=\arrow@iv \advance\angcount@-5 \multiply\angcount@13
   \else
    \ifnum\angcount@<23
     \let\arrowfont@=\arrow@v \advance\angcount@-14 \multiply\angcount@13
    \else
     \let\arrowfont@=\arrow@i \angcount@=117
  \fi\fi\fi
  \ifnum\angcount@=117 \Tcount@=115 \else\Tcount@\angcount@ \fi
 \fi
 \Scount@\Tcount@
 \ifE@
  \ifnum\tcount@=\z@ \advance\Tcount@\tw@ \else\ifnum\tcount@=13
   \advance\Tcount@\tw@ \else \advance\Tcount@\tcount@ \fi\fi
  \ifnum\scount@=\z@ \else \ifnum\scount@=13 \advance\Scount@\thr@@ \else
   \advance\Scount@\scount@ \fi\fi
 \else
  \ifcase\tcount@\advance\Tcount@\thr@@\or\or\advance\Tcount@\thr@@\or
  \advance\Tcount@\tw@\or\advance\Tcount@6 \or\advance\Tcount@7
  \or\or\or\or\or \advance\Tcount@8 \or\advance\Tcount@9 \or
  \advance\Tcount@12 \or\advance\Tcount@\thr@@\fi
  \ifcase\scount@\or\or\advance\Scount@\thr@@\or\advance\Scount@\tw@\or
  \or\or\advance\Scount@4 \or\advance\Scount@5 \or\advance\Scount@10
  \or\advance\Scount@11 \or\or\or\advance\Scount@12 \or\advance
  \Scount@\tw@\fi
 \fi
 \ifcase\arrcount@\or\or\advance\angcount@\@ne\else\fi
 \ifN@ \shifted@=\firsty@ \else\shifted@=-\firsty@ \fi
 \ifE@ \else\advance\shifted@\charht@ \fi
 \goal@=\secondy@ \advance\goal@-\firsty@
 \ifN@\else\multiply\goal@\m@ne\fi
 \setbox\shaft@\hbox{\arrowfont@\char\angcount@}%
 \ifnum\arrcount@=\thr@@
  \getcos@{1.5\p@}%
  \setbox\shaft@\hbox to\wd\shaft@{\arrowfont@
   \rlap{\hskip\dimen@ii
    \smash{\ifNESW@\let\next@\lower\else\let\next@\raise\fi
     \next@\dimen@\hbox{\arrowfont@\char\angcount@}}}%
   \rlap{\hskip-\dimen@ii
    \smash{\ifNESW@\let\next@\raise\else\let\next@\lower\fi
      \next@\dimen@\hbox{\arrowfont@\char\angcount@}}}\hfil}%
 \fi
 \rlap{\smash{\hskip\tocenter@\hskip\firstx@
  \ifnum\arrcount@=\m@ne
  \else
   \ifnum\arrcount@=\thr@@
   \else
    \ifnum\scount@=\m@ne
    \else
     \ifnum\scount@=\z@
     \else
      \setbox\zer@\hbox{\ifnum\angcount@=117 \arrow@v\else\arrowfont@\fi
       \char\Scount@}%
      \ifNESW@
       \ifnum\scount@=\tw@
        \dimen@=\shifted@ \advance\dimen@-\charht@
        \ifN@\hskip-\wd\zer@\fi
        \Nnext@
        \next@\dimen@\copy\zer@
        \ifN@\else\hskip-\wd\zer@\fi
       \else
        \Nnext@
        \ifN@\else\hskip-\wd\zer@\fi
        \next@\shifted@\copy\zer@
        \ifN@\hskip-\wd\zer@\fi
       \fi
       \ifnum\scount@=12
        \advance\shifted@\charht@ \advance\goal@-\charht@
        \ifN@ \hskip\wd\zer@ \else \hskip-\wd\zer@ \fi
       \fi
       \ifnum\scount@=13
        \getcos@{\thr@@\p@}%
        \ifN@ \hskip\dimen@ \else \hskip-\wd\zer@ \hskip-\dimen@ \fi
        \adjust@\shifted@ \advance\adjust@\dimen@ii
        \Nnext@
        \next@\adjust@\copy\zer@
        \ifN@ \hskip-\dimen@ \hskip-\wd\zer@ \else \hskip\dimen@ \fi
       \fi
      \else
       \ifN@\hskip-\wd\zer@\fi
       \ifnum\scount@=\tw@
        \ifN@ \hskip\wd\zer@ \else \hskip-\wd\zer@ \fi
        \dimen@=\shifted@ \advance\dimen@-\charht@
        \Nnext@
        \next@\dimen@\copy\zer@
        \ifN@\hskip-\wd\zer@\fi
       \else
        \Nnext@
        \next@\shifted@\copy\zer@
        \ifN@\else\hskip-\wd\zer@\fi
       \fi
       \ifnum\scount@=12
        \advance\shifted@\charht@ \advance\goal@-\charht@
        \ifN@ \hskip-\wd\zer@ \else \hskip\wd\zer@ \fi
       \fi
       \ifnum\scount@=13
        \getcos@{\thr@@\p@}%
        \ifN@ \hskip-\wd\zer@ \hskip-\dimen@ \else \hskip\dimen@ \fi
        \adjust@\shifted@ \advance\adjust@\dimen@ii
        \Nnext@
        \next@\adjust@\copy\zer@
        \ifN@ \hskip\dimen@ \else \hskip-\dimen@ \hskip-\wd\zer@ \fi
       \fi	
      \fi
  \fi\fi\fi\fi
  \ifnum\arrcount@=\m@ne
  \else
   \loop
    \ifdim\goal@>\charht@
    \ifE@\else\hskip-\charwd@\fi
    \Nnext@
    \next@\shifted@\copy\shaft@
    \ifE@\else\hskip-\charwd@\fi
    \advance\shifted@\charht@ \advance\goal@ -\charht@
    \repeat
   \ifdim\goal@>\z@
    \dimen@=\charht@ \advance\dimen@-\goal@
    \divide\dimen@\tan@i \multiply\dimen@\tan@ii
    \ifE@ \hskip-\dimen@ \else \hskip-\charwd@ \hskip\dimen@ \fi
    \adjust@=\shifted@ \advance\adjust@-\charht@ \advance\adjust@\goal@
    \Nnext@
    \next@\adjust@\copy\shaft@
    \ifE@ \else \hskip-\charwd@ \fi
   \else
    \adjust@=\shifted@ \advance\adjust@-\charht@
   \fi
  \fi
  \ifnum\arrcount@=\m@ne
  \else
   \ifnum\arrcount@=\thr@@
   \else
    \ifnum\tcount@=\m@ne
    \else
     \setbox\zer@
      \hbox{\ifnum\angcount@=117 \arrow@v\else\arrowfont@\fi\char\Tcount@}%
     \ifnum\tcount@=\thr@@
      \advance\adjust@\charht@
      \ifE@\else\ifN@\hskip-\charwd@\else\hskip-\wd\zer@\fi\fi
     \else
      \ifnum\tcount@=12
       \advance\adjust@\charht@
       \ifE@\else\ifN@\hskip-\charwd@\else\hskip-\wd\zer@\fi\fi
      \else
       \ifE@\hskip-\wd\zer@\fi
     \fi\fi
     \Nnext@
     \next@\adjust@\copy\zer@
     \ifnum\tcount@=13
      \hskip-\wd\zer@
      \getcos@{\thr@@\p@}%
      \ifE@\hskip-\dimen@ \else\hskip\dimen@ \fi
      \advance\adjust@-\dimen@ii
      \Nnext@
      \next@\adjust@\box\zer@
     \fi
  \fi\fi\fi}}%
 \iflabel@i
  \rlap{\hskip\tocenter@
  \dimen@\firstx@ \advance\dimen@\secondx@ \divide\dimen@\tw@
  \advance\dimen@\ldimen@i
  \dimen@ii\firsty@ \advance\dimen@ii\secondy@ \divide\dimen@ii\tw@
  \multiply\ldimen@i\tan@i \divide\ldimen@i\tan@ii
  \ifNESW@ \advance\dimen@ii\ldimen@i \else \advance\dimen@ii-\ldimen@i \fi
  \setbox\zer@\hbox{\ifNESW@\else\ifnum\arrcount@=\thr@@\hskip4\p@\else
   \hskip\tw@\p@\fi\fi
   $\m@th\tsize@@\label@i$\ifNESW@\ifnum\arrcount@=\thr@@\hskip4\p@\else
   \hskip\tw@\p@\fi\fi}%
  \ifnum\arrcount@=\m@ne
   \ifNESW@ \advance\dimen@.5\wd\zer@ \advance\dimen@\p@ \else
    \advance\dimen@-.5\wd\zer@ \advance\dimen@-\p@ \fi
   \advance\dimen@ii-.5\ht\zer@
  \else
   \advance\dimen@ii\dp\zer@
   \ifnum\slcount@<6 \advance\dimen@ii\tw@\p@ \fi
  \fi
  \hskip\dimen@
  \ifNESW@ \let\next@\llap \else\let\next@\rlap \fi
  \next@{\smash{\raise\dimen@ii\box\zer@}}}%
 \fi
 \iflabel@ii
  \ifnum\arrcount@=\m@ne
  \else
   \rlap{\hskip\tocenter@
   \dimen@\firstx@ \advance\dimen@\secondx@ \divide\dimen@\tw@
   \ifNESW@ \advance\dimen@\ldimen@ii \else \advance\dimen@-\ldimen@ii \fi
   \dimen@ii\firsty@ \advance\dimen@ii\secondy@ \divide\dimen@ii\tw@
   \multiply\ldimen@ii\tan@i \divide\ldimen@ii\tan@ii
   \advance\dimen@ii\ldimen@ii
   \setbox\zer@\hbox{\ifNESW@\ifnum\arrcount@=\thr@@\hskip4\p@\else
    \hskip\tw@\p@\fi\fi
    $\m@th\tsize@\label@ii$\ifNESW@\else\ifnum\arrcount@=\thr@@\hskip4\p@
    \else\hskip\tw@\p@\fi\fi}%
   \advance\dimen@ii-\ht\zer@
   \ifnum\slcount@<9 \advance\dimen@ii-\thr@@\p@ \fi
   \ifNESW@ \let\next@\rlap \else \let\next@\llap \fi
   \hskip\dimen@\next@{\smash{\raise\dimen@ii\box\zer@}}}%
  \fi
 \fi
}
\def\outnewCD@#1{\def#1{\Err@{\string#1 must not be used within \string\newCD}}}
\newskip\prenewCDskip@
\newskip\postnewCDskip@
\prenewCDskip@\z@
\postnewCDskip@\z@
\def\prenewCDspace#1{\RIfMIfI@
 \onlydmatherr@\prenewCDspace\else\advance\prenewCDskip@#1\relax\fi\else
 \onlydmatherr@\prenewCDspace\fi}
\def\postnewCDspace#1{\RIfMIfI@
 \onlydmatherr@\postnewCDspace\else\advance\postnewCDskip@#1\relax\fi\else
 \onlydmatherr@\postnewCDspace\fi}
\def\predisplayspace#1{\RIfMIfI@
 \onlydmatherr@\predisplayspace\else
 \advance\abovedisplayskip#1\relax
 \advance\abovedisplayshortskip#1\relax\fi
 \else\onlydmatherr@\prenewCDspace\fi}
\def\postdisplayspace#1{\RIfMIfI@
 \onlydmatherr@\postdisplayspace\else
 \advance\belowdisplayskip#1\relax
 \advance\belowdisplayshortskip#1\relax\fi
 \else\onlydmatherr@\postdisplayspace\fi}
\def\PrenewCDSpace#1{\global\prenewCDskip@#1\relax}
\def\PostnewCDSpace#1{\global\postnewCDskip@#1\relax}
\def\newCD#1\endnewCD{%
 \outnewCD@\cgaps\outnewCD@\rgaps\outnewCD@\Cgaps\outnewCD@\Rgaps
 \prenewCD@#1\endnewCD
 \advance\abovedisplayskip\prenewCDskip@
 \advance\abovedisplayshortskip\prenewCDskip@
 \advance\belowdisplayskip\postnewCDskip@
 \advance\belowdisplayshortskip\postnewCDskip@
 \vcenter{\vskip\prenewCDskip@ \Let@ \colcount@\@ne \rowcount@\z@
  \everycr{%
   \noalign{%
    \ifnum\rowcount@=\Rowcount@
    \else
     \global\nointerlineskip
     \getrgap@\rowcount@ \vskip\getdim@
     \global\advance\rowcount@\@ne \global\colcount@\@ne
    \fi}}%
  \tabskip\z@
  \halign{&\global\xoff@\z@ \global\yoff@\z@
   \getcgap@\colcount@ \hskip\getdim@
   \hfil\vrule height10\p@ width\z@ depth\z@
   $\m@th\displaystyle{##}$\hfil
   \global\advance\colcount@\@ne\cr
   #1\crcr}\vskip\postnewCDskip@}%
 \prenewCDskip@\z@\postnewCDskip@\z@
 \def\getcgap@##1{\ifcase##1\or\getdim@\z@\else\getdim@\standardcgap\fi}%
 \def\getrgap@##1{\ifcase##1\getdim@\z@\else\getdim@\standardrgap\fi}%
 \let\Width@\relax\let\Height@\relax\let\Depth@\relax\let\Rowheight@\relax
 \let\Rowdepth@\relax\let\Colwdith@\relax
}
\catcode`\@=\active
\input amsppt.sty
\hsize 12.7cm
\vsize 17.8cm
\magnification=\magstep1
\def\nmb#1#2{#2}         
\def\totoc{}             
\def\idx{}               
\def\ign#1{}             

\redefine\o{\circ}
\define\X{\frak X}
\define\al{\alpha}

\define\de{\delta}

\define\la{\lambda}

\define\si{\sigma}

\define\ph{\varphi}
\define\ch{\chi}
\define\ps{\psi}
\define\om{\omega}

\define\De{\Delta}

\define\La{\Lambda}

\define\Ph{\Phi}
\define\Ps{\Psi}
\define\Om{\Omega}
\predefine\ii\i
\redefine\i{^{-1}}
\define\row#1#2#3{#1_{#2},\ldots,#1_{#3}}
\define\x{\times}
\define\im{\operatorname{im}}
\define\Der{\operatorname{Der}}
\define\Hom{\operatorname{Hom}}
\define\sign{\operatorname{sign}}
\define\Aut{\operatorname{Aut}}
\define\End{\operatorname{End_{\Bbb K}}}
\define\ad{\operatorname{ad}}

\define\Ad{\operatorname{Ad}}
\define\ev{\operatorname{ev}}
\define\Lip{\operatorname{Lip}}
\redefine\L{{\Cal L}}
\def\today{\ifcase\month\or
 January\or February\or March\or April\or May\or June\or
 July\or August\or September\or October\or November\or December\fi
 \space\number\day, \number\year}
\hyphenation{ho-mo-mor-phism}
\topmatter
\title  
The Fr\"olicher-Nijenhuis Bracket  \\
in Non Commutative Differential Geometry \endtitle
\author  Andreas Cap\\
Andreas Kriegl\\
Peter W. Michor  \\
Ji\v r\'\ii{} Van\v zura
\endauthor
\leftheadtext{\smc Cap, Kriegl, Michor, Van\v zura}
\rightheadtext{\smc Non commutative Fr\"olicher-Nijenhuis bracket}
\thanks{Supported by Project P 7724 PHY 
of `Fonds zur F\"orderung der wissenschaftlichen 
Forschung'\hfill}\endthanks
\affil
Institut f\"ur Mathematik, Universit\"at Wien,\\
Strudlhofgasse 4, A-1090 Wien, Austria.\\ \\
Mathematical Institute of the \v CSAV, department Brno,\\ 
Mendelovo n\'am. 1, CS 662 82 Brno, Czechoslovakia
\endaffil
\address{Andreas Cap, Andreas Kriegl, Peter Michor:
Institut f\"ur Mathematik, Universit\"at Wien,
Strudlhofgasse 4, A-1090 Wien, Austria.}\endaddress
\email {michor\@awirap.bitnet} \endemail
\address{Ji\v r\'\ii{} Van\v zura: Mathematical Institute of the \v CSAV,
department Brno, Mendelovo n\'am. 1, CS 662 82 Brno, 
Czechoslovakia}\endaddress
\email {mathmu\@cspuni12.bitnet} \endemail

\endtopmatter

\document

\heading Table of contents \endheading
\noindent Introduction \leaders \hbox to 1em{\hss .\hss }\hfill 
	{\eightrm 1}\par 
\noindent 1. Convenient vector spaces \leaders \hbox to 1em{\hss .\hss }
	\hfill {\eightrm 4}\par 
\noindent 2. Non-commutative differential forms 
	\leaders \hbox to 1em{\hss .\hss }\hfill {\eightrm 7}\par 
\noindent 3. Some related questions 
	\leaders \hbox to 1em{\hss .\hss }\hfill {\eightrm 14}\par 
\noindent 4. The calculus of Fr{\accent "7F o}licher and Nijenhuis 
	\leaders \hbox to 1em{\hss .\hss }\hfill {\eightrm 17}\par 
\noindent 5. Distributions and integrability\leaders 
	\hbox to 1em{\hss .\hss }\hfill {\eightrm 23}\par 
\noindent 6. Bundles and connections 
	\leaders \hbox to 1em{\hss .\hss }\hfill {\eightrm 25}\par 
\noindent 7. Polyderivations and the Schouten-Nijenhuis bracket
	\leaders \hbox to 1em{\hss .\hss }\hfill {\eightrm 29}\par 

\heading\totoc Introduction \endheading

There seems to be an emerging theory of non-commutative differential 
geometry. In the beginning the ideas of non-commutative geometry and 
of non-commutative topology were intended as tools for attacking 
problems in topology, in particular the Novikov conjecture and, more 
generally, the Baum-Connes conjecture. Later on, often motivated by 
physics, one tended to consider `non-commutative spaces' as basic 
structures and to study them in their own right. This is also the 
point of view we adopt in this paper. We carry over to a quite 
general non-commutative setting some of the basic tools of 
differential geometry. From the very beginning we use the setting of 
convenient vector spaces developed by Fr\"olicher and Kriegl. The 
reasons for this are the following: If the non-commutative theory 
should contain some version of differential geometry, a manifold $M$ 
should be represented by the algebra $C^\infty(M,\Bbb R)$ of smooth 
functions on it. The simplest considerations of groups (and quantum 
groups begin to play an important role now) need products, and 
$C^\infty(M\x N,\Bbb R)$ is a certain completion of the algebraic 
tensor product $C^\infty(M,\Bbb R)\otimes C^\infty(N,\Bbb R)$. Now 
the setting of convenient vector spaces offers in its multilinear 
version a monoidally closed category, i.e. there is an appropriate 
tensor product which has all the usual (algebraic) properties with 
respect to bounded multilinear mappings. So multilinear algebra is 
carried into this kind of functional analysis without loss. Moreover 
convenient spaces are the best realm for differentiation which we 
need in section \nmb!{6} to treat a non-commutative version of 
principal bundles.     

We note that all results of this paper also hold in a purely 
algebraic setting: Just equip each vector space with the finest
locally convex topology, then all linear mappings are bounded.
They even remain valid if we take a commutative ring of 
characteristic $\neq 2,3$ instead of the ground field.

In the first section we give a short description of the setting of 
convenient spaces elaborating those aspects which we will need later. 
Then we repeat the usual construction of non-commutative differential 
forms for convenient algebras in the second section.
There we consider triples $(A,\Om_*^A,d)$, where $(\Om^A_*,d)$ is a 
graded differential algebra with $\Om^A_0=A$ and $\Om^A_n=0$ for 
negative $n$. Such a triple is called a \idx{\it quasi resolution} of 
$A$ in the book \cite{Karoubi, 1987}. See in particular 
\cite{Dubois-Violette, 1988} who studies the action of the Lie 
algebra of all derivations on $\Om_*^A$. We will call $(\Om_*^A,d)$ a 
differential algebra for $A$. A universal construction of such an 
algebra $\Om_*^A$ for a commutative algebra $A$ is described in 
\cite{Kunz, 1986}, where it is called the algebra of K\"ahler 
differentials, since apparently this notion was proposed for the 
first time by \cite{K\"ahler, 1953}. The first ones to subsume the 
theory of K\"ahler differentials over a regular affine variety under 
standard homological algebra were \cite{Hochschild, Kostant, 
Rosenberg, 1962}. We present below a non-commutative version of the 
construction of Kunz, since we will need more information. This is 
the construction of \cite{Karoubi, 1982, 1983} which is also used in 
\cite{Connes, 1985}. Connes' contributions started the general 
interest in non-commutative differential geometry. He described the 
Chern character in K-homology coming from Fredholm modules and used 
the universal differential forms as a tool for describing the cyclic 
cohomology of an algebra.

Next we show that the bimodule $\Om_n(A)$ represents the functor of 
the normalized Hochschild $n$-cocyles; this is in principle contained 
in \cite{Connes, 1985}. In the third section we introduce the 
non-commutative version of the Fr\"olicher-Nijenhuis bracket by 
investigating all bounded graded derivations of the algebra of 
differential forms. This bracket is then used to formulate the 
concept of integrability and involutiveness for distributions and to 
indicate a route towards a theorem of Frobenius (the central result 
of usual differential geometry, if there is one). This is then used 
to discuss bundles and connections in the non-commutative setting and 
to go some steps towards a non-commutative Chern-Weil homomorphism. 
In the final section we give a brief description of the 
non-commutative version of the Schouten-Nijenhuis bracket and 
describe Poisson structures.

This work was ignited by a very stimulating talk of Max Karoubi in 
\v Cesky Sternberk in June 1989, and we want to thank him for that.

\newpage
\heading\totoc\nmb0{1}. Convenient vector spaces \endheading

\subheading{\nmb.{1.1}} The traditional differential calculus works 
well for Banach spaces. For more general locally convex spaces a 
whole flock of different theories were developed, each of them rather 
complicated and none really convincing. The main difficulty is that 
the composition of linear mappings stops to be jointly continuous at 
the level of Banach spaces, for any compatible topology. This was the 
original motivation for the development of a whole new field within 
general topology, convergence spaces.

Then in 1982, Alfred Fr\"olicher and Andreas Kriegl presented 
independently the solution to the quest for the right differential 
calculus in infinite dimensions. They joined forces in the further 
development of the theory and the (up to now) final outcome is the 
book \cite{Fr\"olicher, Kriegl, 1988}.

The appropriate spaces for this differential calculus are the 
convenient vector spaces mentioned above. In addition to their
importance for differential calculus these spaces form a category
with very nice properties.

In this section we will sketch the basic definitions and the most 
important results concerning convenient vector spaces and
Fr\"olicher-Kriegl calculus. All locally convex spaces will be 
assumed to be Hausdorff.

\subheading{\nmb.{1.2}. The $c^\infty$-topology} Let $E$ be a 
locally convex vector space. A curve $c:\Bbb R\to E$ is called 
{\it smooth} or $C^\infty$ if all derivatives exist (and are 
continuous) - this is a concept without problems. Let 
$C^\infty(\Bbb R,E)$ be the space of smooth curves. It can be 
shown that $C^\infty(\Bbb R,E)$ does not depend on the locally convex 
topology of $E$, only on its associated bornology (system of bounded 
sets).

The final topologies with respect to the following sets of mappings 
into E coincide:
\roster
\item $C^\infty(\Bbb R,E)$.
\item Lipschitz curves (so that $\{\frac{c(t)-c(s)}{t-s}:t\neq s\}$ 
     is bounded in $E$). 
\item $\{E_B\to E: B\text{ bounded absolutely convex in }E\}$, where 
     $E_B$ is the linear span of $B$ equipped with the Minkowski 
     functional $p_B(x):= \inf\{\la>0:x\in\la B\}$.
\item Mackey-convergent sequences $x_n\to x$ (there exists a sequence 
     $0<\la_n\nearrow\infty$ with $\la_n(x_n-x)$ bounded).
\endroster
This topology is called the $c^\infty$-topology on $E$ and we write 
$c^\infty E$ for the resulting topological space. In general (on the 
space $\Cal D$ of test functions for example) it is finer than the 
given locally convex topology; it is not a vector space topology, 
since addition is no longer jointly continuous. The finest among all 
locally convex topologies on $E$ which are coarser than the 
$c^\infty$-topology is the bornologification of the given locally 
convex topology. If $E$ is a Fr\'echet space, then $c^\infty E = E$.

\subheading{\nmb.{1.3}. Convenient vector spaces} Let $E$ be a 
locally convex vector space. $E$ is said to be a {\it convenient 
vector space} if one of the following equivalent
conditions is satisfied (called $c^\infty$-completeness):
\roster
\item Any Mackey-Cauchy-sequence (so that $(x_n-x_m)$ is Mackey 
     convergent to 0) converges. 
\item If $B$ is bounded closed absolutely convex, then $E_B$ is a 
     Banach space.
\item Any Lipschitz curve in $E$ is locally Riemann integrable.
\item For any $c_1\in C^\infty(\Bbb R,E)$ there is 
     $c_2\in C^\infty(\Bbb R,E)$ with $c_1=c_2'$ (existence of 
     antiderivative).
\endroster

Obviously $c^\infty$-completeness is weaker than 
sequential completeness so any  sequentially complete locally convex
vector space is convenient.
 From \nmb!{1.2}.4 one easily sees that $c^\infty$-closed linear
subspaces of convenient vector spaces are again convenient. We 
always assume that a convenient vector space is equipped with its 
bornological topology.

\proclaim{\nmb.{1.4}. Lemma} Let $E$ be a locally convex space.
Then the following properties are equivalent:
\roster
\item $E$ is $c^\infty$-complete.
\item If $f:\Bbb R\to E$ is scalarwise $\Lip^k$, then $f$ is 
     $\Lip^k$, for $k>1$.
\item If $f:\Bbb R\to E$ is scalarwise $C^\infty$ then $f$ is 
     differentiable at 0.
\item If $f:\Bbb R\to E$ is scalarwise $C^\infty$ then $f$ is 
     $C^\infty$.
\endroster
\endproclaim
Here a mapping $f:\Bbb R\to E$ is called $\Lip^k$ if all partial 
derivatives up to order $k$ exist and are Lipschitz, locally on 
$\Bbb R$. $f$ scalarwise $C^\infty$ means that $\la\o f$ is $C^\infty$  
for all continuous linear functionals on $E$.

This lemma says that on a convenient vector space one can recognize 
smooth curves by investigating compositions with continuous linear 
functionals.

\subheading{\nmb.{1.5}. Smooth mappings} Let $E$ and $F$ be locally 
convex vector spaces. A mapping $f:E\to F$ is called {\it smooth} or 
$C^\infty$, if $f\o c\in C^\infty(\Bbb R,F)$ for all 
$c\in C^\infty(\Bbb R,E)$; so 
$f_*: C^\infty(\Bbb R,E)\to C^\infty(\Bbb R,F)$ makes sense.
Let $C^\infty(E,F)$ denote the space of all smooth mappings from $E$ 
to $F$.

For $E$ and $F$ finite dimensional this gives the usual notion of 
smooth mappings: this has been first proved in \cite{Boman, 1967}.
Constant mappings are smooth. Multilinear mappings are smooth if and 
only if they are bounded. Therefore we denote by $L(E,F)$ the space 
of all bounded linear mappings from $E$ to $F$.

\proclaim{\nmb.{1.6}. Lemma }
For any locally convex space $E$ there is a convenient vector space
$\tilde E$ called the completion of $E$ and a bornological embedding 
$i:E\to \tilde E$, which is characterized by the
property that any bounded linear map from $E$ into an arbitrary 
convenient vector space extends to $\tilde E$.
\endproclaim

\subheading{\nmb.{1.7} }
As we will need it later on we describe the completion in a special
situation:
Let $E$ be a locally convex space with completion $i:E\to \tilde E$,
$f:E\to E$ a bounded projection and $\tilde f:\tilde E\to \tilde E$ the 
prolongation of $i\o f$. Then $\tilde f$ is also a projection and 
$\tilde f(\tilde E)=\ker (Id-\tilde f)$ is a $c^\infty$-closed and
thus convenient linear subspace of $\tilde E$. Using that $f(E)$ is a
direct summand in $E$ one easily shows that $\tilde f(\tilde E)$ is the
completion of $f(E)$. This argument applied to $Id-f$ shows that
$\ker (\tilde f)$ is the completion of $\ker (f)$.

\subheading{\nmb.{1.8}. Structure on $C^\infty(E,F)$} We equip the 
space $C^\infty(\Bbb R,E)$ with the bornologification of the topology 
of uniform convergence on compact sets, in all derivatives 
separately. Then we equip the space $C^\infty(E,F)$ with the 
bornologification of the initial topology with respect to all 
mappings $c^*:C^\infty(E,F)\to C^\infty(\Bbb R,F)$, $c^*(f):=f\o c$, 
for all $c\in C^\infty(\Bbb R,E)$.

\proclaim{\nmb.{1.9}. Lemma } For locally convex spaces $E$ and $F$ 
we have:
\roster
\item If $F$ is convenient, then also $C^\infty(E,F)$ is convenient, 
     for any $E$. The space $L(E,F)$ is a closed linear subspace of 
     $C^\infty(E,F)$, so it is convenient also.
\item If $E$ is convenient, then a curve $c:\Bbb R\to L(E,F)$ is 
     smooth if and only if $t\mapsto c(t)(x)$ is a smooth curve in $F$ 
     for all $x\in E$.
\endroster
\endproclaim

\proclaim{\nmb.{1.10}. Theorem} The category of convenient vector 
spaces and smooth mappings is cartesian closed. So we have a natural 
bijection 
$$C^\infty(E\x F,G)\cong C^\infty(E,C^\infty(F,G)),$$
which is even a diffeomorphism.
\endproclaim

Of course this statement is also true for $c^\infty$-open subsets of 
convenient vector spaces. 

\proclaim{\nmb.{1.11}. Corollary } Let all spaces be convenient vector 
spaces. Then the following canonical mappings are smooth.
$$\align
&\operatorname{ev}: C^\infty(E,F)\x E\to F,\quad 
     \operatorname{ev}(f,x) = f(x).\\
&\operatorname{ins}: E\to C^\infty(F,E\x F),\quad
     \operatorname{ins}(x)(y) = (x,y).\\
&(\quad)^\wedge :C^\infty(E,C^\infty(F,G))\to C^\infty(E\x F,G), 
	\quad \hat f(x,y)=f(x)(y).\\
&(\quad)\spcheck :C^\infty(E\x F,G)\to C^\infty(E,C^\infty(F,G)), 
	\quad \check g(x)(y)=g(x,y).\\
&\operatorname{comp}:C^\infty(F,G)\x C^\infty(E,F)\to C^\infty(E,G)\\
&C^\infty(\quad,\quad):C^\infty(F,F')\x C^\infty(E',E)\to 
     C^\infty(C^\infty(E,F),C^\infty(E',F'))\\
&\qquad (f,g)\mapsto(h\mapsto f\o h\o g)\\
&\prod:\prod C^\infty(E_i,F_i)\to C^\infty(\prod E_i,\prod F_i)
\endalign$$
\endproclaim

\proclaim{\nmb.{1.12}. Theorem} Let $E$ and $F$ be convenient vector 
spaces. Then the differential operator 
$$\gather d: C^\infty(E,F)\to C^\infty(E,L(E,F)), \\
df(x)v:=\lim_{t\to0}\frac{f(x+tv)-f(x)}t,
\endgather$$
exists and is linear and bounded (smooth). Also the chain rule holds: 
$$d(f\o g)(x)v = df(g(x))dg(x)v.$$
\endproclaim

\subheading{\nmb.{1.13} }
The category of convenient vector spaces and bounded linear maps
is complete and cocomplete, so all categorical limits and colimits
can be formed. In particular we can form products and direct sums of
convenient vector spaces.

For convenient vector spaces $E_1$,\dots ,$E_n$ and $F$ we can now
consider the space of all bounded $n$-linear maps, $L(E_1,\dots ,E_n;F)$,
which is a closed linear subspace of $C^\infty(\prod _{i=1}^nE_i,F)$
and thus again convenient. It can be shown that multilinear maps are
bounded if and only if they are partially bounded, i\.e\. bounded in 
each coordinate and that there is a natural isomorphism (of 
convenient vector spaces) $L(E_1,\dots ,E_n;F)\cong
L(E_1,\dots ,E_k;L(E_{k+1},\dots ,E_n;F))$

\proclaim{\nmb.{1.14}. Theorem }
On the category of convenient vector spaces there is a unique tensor product
$\tilde \otimes$ which makes the category symmetric monoidally 
closed, i\.e\. there are natural isomorphisms of convenient vector 
spaces 
$$\alignat2
L(E_1;L(E_2,E_3))&\cong L(E_1\tilde \otimes E_2,E_3),&\quad
E_1\tilde \otimes E_2&\cong E_2\tilde \otimes E_1,\\
E_1\tilde \otimes (E_2\tilde \otimes E_3)&\cong (E_1\tilde \otimes
E_2)\tilde \otimes E_3,&\quad E\tilde \otimes \Bbb R&\cong E.
\endalignat$$
\endproclaim
The tensor product can be constructed as follows: On the algebraic 
tensor product put the finest locally convex topology such that the 
canonical bilinear map from the product into the tensor product is
bounded and then take the completion of this space.

\subheading{\nmb.{1.15}. Remarks } Note that the conclusion of 
theorem \nmb!{1.10} is the starting point of the classical calculus of 
variations, where a smooth curve in a space of functions was assumed 
to be just a smooth function in one variable more.

If one wants theorem \nmb!{1.10} to be true and assumes some other obvious 
properties, then the calculus of smooth functions is already uniquely 
determined.

There are, however, smooth mappings which are not continuous. This is 
unavoidable and not so horrible as it might appear at first sight. 
For example the evaluation $E\x E'\to\Bbb R$ is jointly continuous if 
and only if $E$ is normable, but it is always smooth. Clearly smooth 
mappings are continuous for the $c^\infty$-topology.

For Fr\'echet spaces smoothness in the sense described here coincides 
with the notion $C^\infty_c$ of \cite{Keller, 1974}. This is the 
differential calculus used by \cite{Michor, 1980}, \cite{Milnor, 
1984}, and \cite{Pressley, Segal, 1986}.

\heading\totoc\nmb0{2}. Non-commutative differential forms \endheading

\subheading{\nmb.{2.1}. Axiomatic setting for the algebra of 
differential forms}
Throughout this section we assume that $A$ is a convenient algebra,
i\.e\. $A$ is a convenient vector space together with a bounded 
bilinear associative multiplication $A\times A\to A$. Moreover
we assume that $A$ has a unit $1$.
We consider now a graded associative convenient algebra 
$\Om_*^A=\bigoplus_{p\geq0} \Om_p^A$
where $\Om_0^A=A$ and each $\Om_p^A$ is a convenient vector
space, with a bounded bilinear product 
$:\Om_p^A\x \Om_q^A\to \Om_{p+q}^A$, such that there is a bounded
linear mapping $d=d_p:\Om_p^A\to \Om_{p+1}^A$ with $d^2=0$ and
$d(\om_p \om_q) = d\om_p \om_q + (-1)^p\om_p d\om_q$
for all $\om_p\in\Om_p^A$ and $\om_q\in\Om_q^A$. This mapping is called the 
\idx{\it differential} of $\Om_*^A$. Note that we do not assume that 
the product is graded commutative:
$\om_p \om_q\ne (-1)^{pq}\om_q \om_p$ in general.

Let $\overline{[\Om_*^A,\Om_*^A]_r}$ be the locally convex 
closure of the subspace generated by all graded commutators
$[\om_p,\om_q]:= \om_p \om_q -(-1)^{pq}\om_q \om_p$ with
$p+q=r$. We put $\bar \Om_r^A:= 
\Om_r^A/\overline{[\Om_*^A,\Om_*^A]_r}$ and 
we let $T:\Om_r^A\to \bar\Om_r^A$ be the projection which will be 
called the \idx{\it graded trace} of $\Om_*^A$.

Since we have $d([\om_p,\om_q])=[d\om_p,\om_q]+(-1)^p[\om_p,d\om_q]$, 
the differential passes to $\bar\Om_*^A$ and still satisfies 
$d^2=0$. The separated homology of this quotient complex is called the 
\idx{\it non-commutative De Rham homology} of $\Om_*^A$ or of $A$, 
if $\Om_*^A$ is clear. We denote it by
$$H\bar\Om_p^A=\bar H_p^A
     =\ker (d:\bar\Om_p^A\to\bar\Om_{p+1}^A)/
     \overline{\im (d:\bar\Om_{p-1}^A\to\bar\Om_p^A)}.$$

\subheading{\nmb.{2.2}. Derivations} Let $M$ be a convenient bimodule
over the convenient algebra $A$, i\.e\. $M$ is a convenient vector 
space together with two bounded homomorphisms of unital algebras
$\la :A\to L(M,M)$ and $\rho :A^{op}\to L(M,M)$, where $A^{op}$ 
denotes the opposite algebra to $A$, such that for $a,b\in A$ we have 
$\la (a)\o\rho (b)=\rho (b)\o\la (a)$. We will write $am$ for $\la (a)(m)$
and $ma$ for $\rho (a)(m)$. This definition is equivalent to having
bounded bilinear maps $\la :A\times M\to M$ and $\rho :M\times A\to M$,
which satisfy the usual axioms. A (bounded) \idx{\it derivation} of $A$
in $M$ is a bounded linear mapping $D:A\to M$ such that 
$D(ab)=D(a)b+aD(b)$ for all $a$, $b\in A$.
We denote by $\Der(A;M)$ the vector space of all derivations 
of $A$ into $M$. This is obviously a closed linear subspace of
$L(A,M)$ and thus a convenient vector space. If $A$ is commutative, then
$\Der(A;M)$ is again an $A$-module. 

The vector space $\Der(A;A)$ is a convenient Lie algebra where the
bracket is the commutator. It is an $A$-module if and only if $A$ is 
commutative.

\subheading\nofrills{\bf\nmb.{2.3}. The algebra of dual numbers}
{} of a convenient algebra $A$ with respect to a convenient
$A$-bimodule $M$ is the semidirect product $A\circledS M$, 
i.e. the convenient vector space $A\x M$ with the bounded bilinear 
multiplication 
$(a_1,m_1)(a_2,m_2):= (a_1a_2,a_1m_2+m_1a_2)$.
This is an associative convenient algebra with unit $(1,0)$.

\proclaim{\nmb.{2.4}. Lemma} The bounded derivations from $A$ into the 
$A$-bimodule $M$ correspond exactly to the bounded algebra 
homomorphisms $\ph:A\to A\circledS M$ satisfying 
$pr_1\o\ph=Id_A$. \qed
\endproclaim

\subheading{\nmb.{2.5}. Universal derivations} A bounded derivation $D:A\to M$
into a bimodule $M$ is 
called \idx{\it universal\ign{ derivation}} if the following holds:
\roster
\item"" For any bounded derivation $D':A\to N$ into a convenient $A$-bimodule
     $N$ there is a unique bounded $A$-bimodule homomorphism $\Ph:M\to N$
     such that $D'=\Ph\o D$.
\endroster
Of course for any two universal derivations $D_1:A\to M_1$ and 
$D_2:A\to M_2$ there is a unique $A$-bimodule isomorphism 
$\Ph:M_1\to M_2$ such that $D_2=\Ph\o D_1$. So a universal derivation 
is unique up to canonical isomorphism.

\proclaim{Lemma} For every convenient algebra $A$ there exists a universal 
derivation which we denote by $d:A\to \Om_1(A)$.
\endproclaim
\demo{Proof}
First we define an $A$-bimodule structure on $A\tilde \otimes A$ as 
follows: Let $(a,b)\mapsto a\otimes b:A\times A\to A\tilde \otimes A$ be
the canonical bilinear map. Now consider the map 
$\bar \la :A\to L(A\times A,A\tilde \otimes A)$ defined by
$\bar\la (a)(b,c):=ab\otimes c$. Obviously the map $\bar\la$ has values in
the space $L(A,A;A\tilde \otimes A)$ of bilinear maps and thus we can 
compose it with the isomorphisms of \nmb!{1.13} and \nmb!{1.14} to get
$\la :A\to L(A\tilde \otimes A,A\tilde \otimes A)$ which is easily 
seen to be an algebra homomorphism.
Similarly we define $\rho :A\to L(A\tilde \otimes A,A\tilde \otimes A)$
using $\bar \rho (a)(b,c):=b\otimes ca$.

The multiplication on $A$ induces a bounded linear map 
$\mu :A\tilde \otimes A\to A$ which is an $A$-bimodule homomorphism
by associativity. Thus $\Om _1(A):=\ker (\mu )$ is a convenient 
$A$-bimodule.

Next we define $d:A\to \Om _1(A)$ by $d(a):=1\otimes a-a\otimes 1$. Obviously
$d$ is a bounded derivation.

To see that this derivation is universal let $D:A\to M$ be a bounded
derivation from $A$ into a convenient $A$-bimodule $M$. 
Let $\bar\Phi :A\times A\to M$ be the map defined by $\bar\Phi (a,b):=aD(b)$.
Then $\bar\Phi$ is obviously bilinear and bounded and thus it induces a
bounded linear map $\Phi :A\tilde \otimes A\to M$, whose restriction to
$\Om _1(A)$ we also denote by $\Phi$.
As any derivation vanishes on $1$ we get:
$$(\Phi \o d)(a)=\Phi (1\otimes a-a\otimes 1)=1D(a)-aD(1)=D(a)$$
So it remains to show that $\Phi$ is a bimodule homomorphism.
For $a,b,c\in A$ we get:
$(\Phi \o \la (a))(b\otimes c)=\Phi (ab\otimes c)=abD(c)=a(\Phi (b\otimes c))$ 
and thus $\Phi:A\tilde \otimes A\to M$ is a homomorphism of left modules.

On the other hand $(\Phi \o \rho (a))(b\otimes c)=bD(ca)=(bD(c))a+bcD(a)$
and thus we get the identity 
$(\Phi \o \rho (a))(x)=(\Phi (x))a+\mu (x)D(a)$ for all 
$x\in A\tilde \otimes A$ and so 
$\Phi :\Om _1(A)\to M$ is a homomorphism of right modules, too.
\qed\enddemo

\proclaim{\nmb.{2.6}. Corollary} For an $A$-bimodule $M$ the canonical 
linear mapping 
$$\gather
d^*:\Hom^A_A(\Om_1(A),M)\to \Der(A;M)\\
\ph\mapsto \ph\o d
\endgather$$
is an isomorphism of convenient vector spaces, where $\Der(A;M)$ 
carries the structure described in \nmb!{2.2}, while the space 
$\Hom^A_A(\Om_1(A),M)$ of 
bounded bimodule homomorphisms is considered as a 
closed linear subspace of $L(\Om_1(A),M)$.
In particular we have 
$\Hom^A_A(\Om_1(A),A)\cong \Der(A;A)$.
\endproclaim

\demo{Proof}
Since $d$ is bounded and linear so is $d^*$. In the proof of the 
lemma above we saw that the inverse to $d^*$ is given by mapping
$D$ to the prolongation of $\ell \o (Id\times D)$, where $\ell$
denotes the left action of $A$ on $M$ and this map is easily seen
to be bounded.
\qed\enddemo

\proclaim{\nmb.{2.7}. Lemma }
Let $A$ be a convenient algebra, $M$ a convenient right $A$-module
and $N$ a convenient left $A$-module.
\roster
\item There is a convenient vector space $M\tilde \otimes _AN$ and a
     bounded bilinear map $b:M\times N\to M\tilde \otimes _AN$,
     $(m,n)\mapsto m\otimes _An$ such that
     $b(ma,n)=b(m,an)$ for all $a\in A$, $m\in M$ and $n\in N$ which has
     the following universal property: If $E$ is a convenient vector space
     and $f:M\times N\to E$ is a bounded bilinear map such that 
     $f(ma,n)=f(m,an)$ then there is a unique bounded linear map
     $\tilde f:M\tilde \otimes _AN\to E$ with $\tilde f\o b=f$.
\item Let $L^A(M,N;E)$ denote the space of all bilinear bounded maps
     $f:M\times N\to E$ having the above property, which is a closed
     linear subspace of $L(M,N;E)$. Then we have an isomorphism of 
     convenient vector spaces $L^A(M,N;E)\cong L(M\tilde \otimes _AN,E)$.
\item If $B$ is another convenient algebra such that $N$ is a 
     convenient right $B$-module and such that the actions of $A$
     and $B$ on $N$ commute, then $M\tilde \otimes _AN$ is in a canonical
     way a convenient right $B$-module.
\item If in addition $P$ is a convenient left $B$-module then there
     is a natural isomorphism of convenient vector spaces
$$M\tilde \otimes _A(N\tilde \otimes _BP)\cong (M\tilde \otimes 
     _AN)\tilde \otimes _BP$$
\endroster
\endproclaim
\demo{Proof}
We construct $M\tilde \otimes _AN$ as follows: Let $M\otimes N$ be 
the algebraic tensor product of $M$ and $N$ equipped with the (bornological)
topology mentioned in \nmb!{1.14} and let $V$ be the locally 
convex closure of the subspace generated by all elements of the form
$ma\otimes n-m\otimes an$ and define $M\tilde \otimes _AN$ to be the
completion of $M\otimes _AN:=(M\otimes N)/V$. As $M\otimes N$ has the
universal property that bounded  bilinear maps from $M\times N$ into
arbitrary locally convex spaces induce bounded and hence continuous
linear maps on $M\otimes N$, \therosteritem{1} is clear.

\therosteritem{2}:
By \therosteritem {1} the bounded linear map 
$b^*:L(M\tilde \otimes _AN,E)\to L^A(M,N;E)$ is a bijection. Thus it
suffices to show that its inverse is bounded, too. From \nmb!{1.14}
we get a bounded linear map $\ph :L(M,N;E)\to L(M\otimes N,E)$ which is 
inverse to the map induced by the canonical bilinear map.
Now let $L^{\text{ann }V}(M\otimes N,E)$ be the closed linear subspace of 
$L(M\otimes N,E)$ consisting of all maps which annihilate $V$.
Restricting $\ph$ to $L^A(M,N;E)$ we get a bounded linear map
$\ph :L^A(M,N;E)\to L^{\text{ann }V}(M\otimes N,E)$.

Let $\ps :M\otimes N\to M\otimes _AN\to M\tilde \otimes _AN$ be the
composition of the canonical projection with the inclusion into the
completion. Then $\ps$ induces a well defined linear map
$\hat \ps :L^{\text{ann }V}(M\otimes N,E)\to L(M\tilde \otimes _AN,E)$ and
$\hat \ps \o \ph$ is inverse to $b^*$. So it suffices to show that 
$\hat \ps$ is bounded.

This is the case if and only if the associated map
$L^{\text{ann }V}(M\otimes N,E)\times (M\tilde \otimes _AN)\to E$ is bounded.
This in turn is equivalent to boundedness of the associated map
$M\tilde \otimes _AN\to L(L^{\text{ann }V}(M\otimes N,E),E)$.
But this is just the prolongation to the completion of the map
$M\otimes _AN\to L(L^{\text{ann }V}(M\otimes N,E),E)$ which sends $x$ to the 
evaluation at $x$ and this map is clearly bounded.

\therosteritem{3}:
Let $\rho :B^{op}\to L(N,N)$ be the right action of $B$ on $N$ and 
let $\Phi :L^A(M\times N,M\tilde \otimes _AN)\cong
L(M\tilde \otimes _AN,M\tilde \otimes _AN)$ be the isomorphism 
constructed in \therosteritem{2}. We define the right module structure
on $M\tilde \otimes _AN$ as:
$$\multline
B^{op} @>\rho >> L(N,N) @>Id\times .>> L(M\times N,M\times N) @>b_*>>\\
@>>> L^A(M,N;M\tilde \otimes _AN) @>\Phi >>
L(M\tilde \otimes _AN,M\tilde \otimes _AN)
\endmultline$$
This map is obviously bounded and easily seen to be an algebra 
ho\-mo\-mor\-phism.

\therosteritem{4}:
Straightforward computations show that both spaces have the following
universal property: For a convenient vector space $E$ and a trilinear
map $f:M\times N\times P\to E$ which satisfies $f(ma,n,p)=f(m,an,p)$
and $f(m,nb,p)=f(m,n,bp)$ there is a unique linear map prolonging $f$.
\qed\enddemo

\subheading{\nmb.{2.8}. Homomorphisms of differential algebras}
Let $\ph :A\to B$ be a homomorphism of convenient algebras, let 
$(\Om^A,d^A)$ be a differential algebra for $A$ in the sense of 
\nmb!{2.1}, and let $(\Om^B,d^B)$ be one for $B$. 

By a $\ph$-homomorphism $\Phi:\Om^A\to\Om^B$ we mean a bounded
homomorphism of graded differential algebras such that 
$\Phi_0=\ph :\Om^A_0=A\to B=\Om^B_0$.

\proclaim{\nmb.{2.9}. Theorem. Existence of the universal graded 
differential algebra}
For each convenient algebra $A$ there is a convenient graded differential
algebra $(\Om(A),d)$ for $A$ with the following property:
\roster
\item"" For any bounded homomorphism $\ph:A\to B$ of convenient algebras and 
for any convenient graded differential algebra $(\Om^B,d^B)$ for $B$ there
exists a unique $\ph$-homomorphism $\Om(A)\to \Om^B$.
\endroster
\endproclaim

\demo{Proof}
Put $\Om _0(A)=A$ and $\Om _k(A):=\Om _1(A)\tilde \otimes _A\dots
\tilde \otimes _A\Om _1(A)$ ($k$ factors). Then each $\Om _k(A)$ is 
a convenient $A$-bimodule by \nmb!{2.7}.3, which also defines 
the multiplication with elements of $\Om _0(A)$. For $k,\ell >0$
we define the multiplication as the canonical bilinear map
$$\Om _k(A)\times \Om _{\ell}(A)\to \Om _k(A)\tilde \otimes _A\Om
_{\ell}(A)\cong \Om _{k+\ell}(A)$$
Thus $\Om (A)=\bigoplus _k\Om _k(A)$ is a convenient graded algebra.

{\bf Claim.} There is an isomorphism $\Om_1(A)\cong A\tilde \otimes
(A/\Bbb R)$ of convenient vector spaces. \newline
Consider the embedding $i:\Bbb R\to A$ and the projection 
$p:A\to A/\Bbb R$, denoted also by $p(a)=:\bar a$. 
We consider the following diagram, where the 
horizontal and the vertical sequences are exact:
$$\cgaps{.5;1;1;.5}\rgaps{.5;1;1;.5}\CD
  @.   @.  0         @.    @. \\
@.  @.     @AAA      @.    @.\\
  @. A @=  A         @.    @.  \\
@.     @|                              @A\mu AA @. @.\\
0 @>>> A\tilde \otimes \Bbb R @>Id\tilde \otimes  i>> A\tilde \otimes
A @>Id\tilde \otimes  p>> A\tilde \otimes  (A/\Bbb R) @>>> 0 \\
@.  @.  @AAA  @.  @.  \\
  @.  @.  \Om_1(A)  @.  @.\\
@.  @.  @AAA  @.  @. \\
  @.  @.  0 @.  @. 
\endCD$$
The vertical sequence is splitting: $a\mapsto a\otimes 1$ is a 
section for $\mu$ and the prolongation of $(a,b)\mapsto a\,d(b)$ is a
retraction onto $\Om_1(A)$ which even factors over $Id\tilde \otimes p$,
since by \nmb!{1.7} the space $\Om _1(A)$ is the completion of the 
kernel of the prolongation of the multiplication map to $A\otimes A$. So
we may invert all arrows of the vertical sequence and the two sequences
are isomorphic as required.

{\bf Claim.} There is an isomorphism of convenient vector spaces
$$A\tilde \otimes \oversetbrace \text{$k$-times} 
\to{A/{\Bbb R}\tilde \otimes \cdots \tilde \otimes A/\Bbb R} \to \Om_k(A)$$
which is induced by the map $(a_0,\bar a_1,\dots ,\bar a_k)\mapsto
a_0da_1\otimes _Ada_2\otimes _A\dots \otimes _Ada_k$.
This is a direct consequence of the last claim and lemma \nmb!{2.7}.

We now define $d:\Om_k(A)\to \Om_{k+1}(A)$ by 
$d(a)=1\otimes a-a\otimes 1$ for $a\in \Om_0(A)= A$ and for $k>0$ as the 
mapping defined on $\Om _k(A)\cong A\tilde \otimes A/\Bbb R\tilde \otimes
\dots \tilde \otimes A/\Bbb R$ which is associated to:
$$\gather
(a_0,\bar a_1,\dots ,\bar a_k)\mapsto 1\otimes \bar a_0\otimes \bar
a_1\otimes \dots \otimes \bar a_k\\
A\times (A/\Bbb R)^k\to A\tilde \otimes A/\Bbb R\tilde \otimes \dots
\tilde \otimes A/\Bbb R\cong \Om _{k+1}(A)
\endgather$$
Let us show now that $d$ is a graded derivation: We have to show that
for $\om _k\in \Om _k(A)$ and $\om _{\ell}\in \Om _{\ell}(A)$ we have
$d(\om _k\om _{\ell})=d(\om _k)\om _{\ell}+(-1)^k\om _kd(\om _{\ell})$.
We proceed by induction on $k$. By the claim above it suffices to 
check the identity for elements of $A\times (A/\Bbb R)^i$. For $k=0$
we have 
$a(b_0,\bar b_1,\dots ,\bar b_{\ell})=(ab_0,\bar b_1,\dots ,\bar b_{\ell})$
which is mapped by $d$ to the element
$(1,\overline{ab_0},\bar b_1,\dots ,\bar b_{\ell})$ which under the isomorphism
with $\Om _{\ell}(A)$ goes to
$d(ab_0)\otimes _Adb_1\otimes _A\dots \otimes _Adb_{\ell}$ so the result 
follows from the derivation property of $d:A\to \Om _1(A)$.

In the general case we first see that using this derivation property
again, the product of
$(a_0,\bar a_1,\dots ,\bar a_k)$ and $(b_0,\bar b_1,\dots ,\bar b_{\ell})$
in $\Om _{k+\ell}(A)$ can be written as
$$
\multline
a_0da_1\otimes _A \dots \otimes _A da_{k-1}\otimes _A
d(a_kb_0)\otimes _A db_1\otimes _A \dots \otimes _A db_{\ell}-\\
-(a_0da_1\otimes _A \dots \otimes _A da_{k-1})(a_kdb_0\otimes _A
db_1\otimes _A \dots \otimes _A db_{\ell})
\endmultline
$$
and from this the result follows easily using the induction 
hypothesis.

So let us turn to the universal property. Let $B$ be a convenient
algebra, $(\Om ^B,d^B)$ a convenient differential algebra for $B$
and $\ph :A\to B$ a bounded homomorphism of algebras.
Via $\ph$ and the multiplication of $\Om ^B$ all spaces $\Om _i^B$
are convenient $A$-bimodules. 

As $d^B$ is a graded derivation the map $d^B\o \ph :A\to \Om_1^B$
is a derivation. Thus by the universal property of $\Om _1(A)$
we get a unique bounded bimodule homomorphism $\ph _1:\Om _1(A)\to
\Om _1^B$. Thus for $a\in A$ and $\om \in \Om _1(A)$ we have
$\ph _1(a\om )=\ph (a)\ph _1(\om )$ and 
$\ph _1(\om a)=\ph _1(\om )\ph (a)$.
Consider the map $f:(\Om _1(A))^k\to \Om _k^B$ defined by
$f(\om _1,\om _2,\dots ,\om _k):=\ph _1(\om _1)\ph _1(\om _2)\dots
\ph _1(\om _k)$ which is obviously bounded and $k$-linear.
Moreover as $\ph _1$ is a bimodule homomorphism we get
$f(\dots ,\om _ia,\om _{i+1},\dots )=f(\dots ,\om _i,a\om _{i+1},\dots )$.
Thus there is a unique prolongation of $f$ to $\Om _k(A)$ which
we define to be $\ph _k$. From this definition it is obvious that
the maps $\ph _i$ form a bounded homomorphism of graded algebras.

The composition:
$$\align
A\times &A/\Bbb R\times \dots \times A/\Bbb R\to A\tilde \otimes A/\Bbb R
\dots \tilde \otimes A/\Bbb R\cong \Om _k(A) @>\ph _k>> \Om _k^B\\
\intertext{is given by} 
&(a_0,\bar a_1,\dots ,\bar a_k)\mapsto a_0da_1\otimes _Ada_2\otimes _A\dots
\otimes _Ada_k\mapsto\\
&\ph (a_0)\ph _1(da_1)\dots \ph _1(da_k)=\ph (a_0)d^B(\ph (a_1))
\dots d^B(\ph (a_k))
\endalign$$
and this element is mapped by $d^B$ to $d^B(\ph (a_0))d^B(\ph (a_1))
\dots d^B(\ph (a_k))$. This shows that 
$\ph _{k+1}\o d=d^B\o \ph _k$
\qed\enddemo

\proclaim{\nmb.{2.10}. Corollary}
The construction $A \mapsto \Om_*(A)$ defines a covariant functor 
from the category of convenient algebras with unit to the category of 
convenient graded differential algebras.
\endproclaim
So for a bounded algebra homomorphism $f:A\to B$ we denote by 
$\Om_*(f):\Om_*(A)\to \Om_*(B)$ its universal prolongation.

\heading\totoc\nmb0{3}. Some related questions \endheading

In the following we treat two questions which arise naturally in the 
context of section \nmb!{2} but which are not relevant for the 
developments afterwards.

\subheading{\nmb.{3.1}} The kernel of the multiplication 
$\mu: A\tilde\otimes A\to A$ is the very important space $\Om_1(A)$. 
What about the analogue with more factors?

\proclaim{Proposition} Let $A$ be a convenient algebra 
with unit. Then the kernel of the $n$-ary multiplication 
$\mu^n:A^{\tilde\otimes n}\to A$ is the subspace
$$\sum_{i=0}^{n-2}A^{\tilde\otimes i}\tilde\otimes\Om_1(A)
	\tilde\otimes A^{\tilde\otimes (n-2-i)}
	\subset A^{\tilde\otimes n}.$$
\endproclaim

\demo{Proof}
Note that $\mu^2=\mu:A\tilde\otimes A\to A$. We prove the assertion by 
induction on $n$. Consider the following commutative diagram:
$$\cgaps{0.3;0.5;0.8;0.3}\rgaps{0.7;1;1;0.7}\newCD
 & & 0 @(0,-1) & 0 @(0,-1) & \\
 & & A^{\tilde\otimes(n-1)}\tilde\otimes \Om_1(A) @(0,-1) @(1,0) 
	& \Om_1(A) @(1,0) @(0,-1) & 0	\\
0 @(1,0) & \left(\sum_{i=0}^{n-2}A^{\tilde\otimes i}\tilde\otimes \Om_1(A) 
	\tilde\otimes A^{\tilde\otimes (n-2-i)} \right)\tilde\otimes A	@(1,0) @(0,-1)
	& A^{\tilde\otimes(n+1)} 
	@()\L{A^{\tilde\otimes(n-1)}\tilde\otimes \mu}@(0,-1) 
     @()\L{\mu^n\tilde\otimes A}@(1,0) @()\L{\mu^{n+1}}@(1,-1)
	& A\tilde\otimes A @(1,0) @()\l{\mu}@(0,-1) & 0 \\
0 @(1,0) & \sum_{i=0}^{n-2}A^{\tilde\otimes i}\tilde\otimes \Om_1(A) 
	\tilde\otimes A^{\tilde\otimes (n-2-i)} @(0,-1) @(1,0) 
	& A^{\tilde\otimes n} @()\L{\mu^n}@(1,0) @(0,-1) & A @(1,0) @(0,-1) & 0 \\
 & 0 & 0 & 0 &
\endnewCD$$
The right hand column is the defining sequence for $\Om_1(A)$ and it 
is splitting. The middle column being the right hand one tensored 
with $A^{\tilde\otimes(n-1)}$ from the left is then again splitting 
and thus exact. The bottom row is exact by the induction hypothesis 
and is also splitting since $\mu^n$ admits many obvious sections. The 
middle row is the bottom one tensored with $A$ from the right and it 
is again splitting and thus exact. The left hand side vertical arrow 
is multiplication from the right. The top horizontal arrow is total 
multiplication onto the left of $\Om_1(A)$.

Let us now take an element $x\in A^{\tilde\otimes (n+1)}$ which is 
in the kernel of $\mu^{n+1}$.
Then a simple diagram chasing shows that $x$ is in 
the sum of the two subspaces of $A^{\tilde\otimes (n+1)}$ which are 
above and to the left.
The converse is trivial, so the result follows.
\qed\enddemo

\subheading{\nmb.{3.2}} We have seen in \nmb!{2.6} that $\Om_1(A)$ is 
the representing object for the functor $\Der(A,\quad)$ on the 
category of $A$-bimodules. Which functor is represented by 
$\Om_n(A)$?

Recall that 
$\Om_n(A) = \Om_1(A)\tilde\otimes_A\dots\tilde\otimes_A\Om_1(A)$ 
($n$ times). We consider the $n$-linear 
mapping 
$$\gather
d^n: A^n\to (A/\Bbb R)^n\to \Om_n(A),\\
d^n(a_1,\dots a_n) := da_1\otimes_A \dots \otimes_A da_n.
\endgather$$ 
We view it as a Hochschild cochain which is bounded as a multilinear 
mapping and normalized, i\. e\. it factors to $(A/\Bbb R)^n$.
It is well known that the normalized Hochschild complex leads to the 
usual Hochschild cohomology, see \cite{Cartan, Eilenberg, 1956, p\. 
176}.

\proclaim{Lemma} The mapping $d^n$ is a normalized and bounded 
Hochschild cocycle with values in the $A$-bimodule $\Om_n(A)$.
\endproclaim

\demo{Proof}
By definition of the right $A$-module structure on $\Om_n(A)$ we have 
$$\align
d^n(a_1&,\dots,a_n)a_{n+1} 
= (da_1\otimes_A \dots \otimes_A da_n)a_{n+1} \\
&= da_1\otimes_A \dots \otimes_A d(a_na_{n+1}) - 
(da_1\otimes_A \dots \otimes_A da_{n-1})a_n\otimes_A da_{n+1}\\
&= d^n(a_1,\dots,a_na_{n+1}) - d^n(a_1,\dots,a_{n-1}a_n,a_{n+1})\\
&\quad + (da_1\otimes_A\dots\otimes_A da_{n-2})a_{n-1}\otimes_A 
	(da_n\otimes_A da_{n+1})\\
&= \ldots\\
&= \sum_{i=1}^n(-1)^{n-i} d^n(a_1,\dots,a_ia_{i+1},\dots,a_{n+1})
+ (-1)^n a_1 d^n(a_2,\dots,a_{n+1}),
\endalign$$
and thus as required
$$\align
0 &= a_1 d^n(a_2,\dots,a_{n+1}) 
	+ \sum_{i=1}^n(-1)^i d^n(a_1,\dots,a_ia_{i+1},\dots,a_{n+1}) \\
&\quad +(-1)^{n+1}d^n(a_1,\dots,a_n)a_{n+1} \\
&=: (\de d^n)(a_1,\dots,a_{n+1}),
\endalign$$
where $\de$ denotes the usual Hochschild coboundary operator.
\qed\enddemo

\proclaim{\nmb.{3.3}. Proposition} Let $M$ be an $A$-bimodule.  Then 
the mapping 
$$(d^n)^*: \Hom_A^A(\Om_n(A),M)\to \bar Z^n(A, M)$$ 
is an isomorphism onto the space of all normalized and bounded 
Hochschild cocycles with values in $M$.
\endproclaim

\demo{Proof}
Clearly for any bimodule homomorphism $\Ph:\Om_n(A)\to M$ the 
$n$-linear mapping $\Ph\o d^n:\bar A^n\to M$ is a normalized and 
bounded Hochschild cocycle.
Let us assume conversely that $c:A^n\to M$ is a normalized bounded cocycle. 
In the proof of \nmb!{2.9} we got a natural isomorphism of 
convenient vector spaces
$$A\tilde \otimes \oversetbrace \text{$k$-times} 
\to{A/{\Bbb R}\tilde \otimes \cdots \tilde \otimes A/\Bbb R} \to \Om_k(A)$$
which is given by 
$a_0 \otimes\bar a_1\otimes\dots \otimes\bar a_k\mapsto
a_0da_1\otimes _Ada_2\otimes _A\dots \otimes _Ada_k$.
Using this we define 
$\Ph_c:\Om_n(A)\to M$ by 
$\Ph_c(a_0da_1\dots da_n):= a_0c(a_1,\dots,a_n)$.
Then clearly $\Ph\o d^n=c$. Obviously $\Ph_c$ is a homomorphism of 
left $A$-modules and from the definition of the right $A$-module 
structure on $\Om_n(A)$ we see that $\de c=0$ translates into 
$\Ph_c$ being a right module homomorphism, by a computation which is 
completely analogous to the one in the proof of \nmb!{3.2}.
Obviously both constructions are bounded.
\qed\enddemo

\subheading{\nmb.{3.4}} Is it possible to recognize the Hochschild 
coboundaries in the description 
$\bar Z^n(A, M)\cong\Hom_A^A(\Om_n(A),M)$?

In order to answer this question we consider the canonical normalized 
mapping, where $a\mapsto \bar a$ is the quotient mapping 
$A\to A/\Bbb R$:
$$\gather
\ph: A^{n-1}\to A\tilde\otimes\oversetbrace{n-1}\to
{(A/\Bbb R)\tilde\otimes\dots\tilde\otimes(A/\Bbb R)}\tilde\otimes A\\
\ph(a_1,\dots,a_{n-1}):= 
1\otimes\bar a_1\otimes\dots\otimes\bar a_{n-1}\otimes 1
\endgather$$
Then $\partial \ph\in 
\bar Z^n(A;A\tilde\otimes(A/\Bbb R)^{\tilde\otimes(n-1)}\tilde\otimes A)$ 
is given by 
$$\align
\partial\ph(a_1,\dots,a_n) 
&= a_1 \ph(a_2,\dots,a_n) 
	+ \sum_{i=1}^{n-1}(-1)^i \ph(a_1,\dots,a_ia_{i+1},\dots,a_n) \\
&\qquad +(-1)^n\ph(a_1,\dots,a_{n-1})a_n\\
&= a_1\otimes\bar a_2\otimes\dots\otimes\bar a_n\otimes 1\\
&\qquad	+ \sum_{i=1}^{n-1}(-1)^i 
    	1\otimes\bar a_1\otimes\dots\otimes\overline{a_ia_{i+1}}\otimes
	\dots\otimes\bar a_n\otimes 1\\
&\qquad +(-1)^n 1\otimes\bar a_1\otimes\dots\otimes\bar a_{n-1}\otimes a_n.
\endalign$$
By proposition \nmb!{3.3} there exists a unique bimodule homomorphism 
$I:\Om_n(A)\to A\tilde\otimes(A/\Bbb R)^{\tilde\otimes(n-1)}\tilde\otimes A$
such that $\partial \ph = I\o d^n$.

A short computation (again essentially the same as in the proof of 
lemma \nmb!{3.2}) shows that this bimodule homomorphism $I$ coincides 
with the following composition of canonical mappings:
$$\multline
\Om_n(A) = \Om_1(A)\tilde\otimes_A\dots\tilde\otimes_A\Om_1(A)
	\East{i\otimes\dots\otimes i}{} \\
\East{i\otimes\dots\otimes i}{} (A\tilde\otimes A)\tilde\otimes_A
	\dots\tilde\otimes_A(A\tilde\otimes A)
\cong \oversetbrace{n+1}\to {A\tilde\otimes\dots\tilde\otimes A}
	\to A\tilde\otimes \oversetbrace{n-1}\to 
	{(A/\Bbb R)\tilde\otimes\dots\tilde\otimes(A/\Bbb R)}\tilde\otimes A,
\endmultline$$
where $i$ is the injection $\Om_1(A)= \ker \mu \to A\tilde\otimes A$.

\proclaim{\nmb.{3.5}. Proposition} Let $\Ph:\Om_n(A)\to M$ be a 
bimodule homomorphism. Then the corresponding normalized Hochschild cocycle 
$\Ph\o d^n$ is a coboundary if and only if $\Ph$	factors 
over $I$ to a bimodule homomorphism
$\tilde\Ph:A\tilde\otimes(A/\Bbb R)^{\tilde\otimes(n-1)}\tilde\otimes A 
\to M$, 
so that $\Ph = \tilde\Ph\o I$. 

In more details: for any bimodule homomorphism 
$\Ps:A\tilde\otimes(A/\Bbb R)^{\tilde\otimes(n-1)}\tilde\otimes A \to M$ 
we have $\Ps\o I\o d^n = \partial \ps$ where the normalized bounded 
cochain
$\ps:A^{n-1}\to M$ 
is given by 
$$\ps(a_1,\dots, a_{n-1}) = 
\Ps(1\otimes\bar a_1\otimes\dots\otimes\bar a_{n-1}\otimes 1).$$
\endproclaim

\demo{Proof}
Let $\Ph\o d^n$ be a coboundary. Then there is an $(n-1)$-linear 
mapping $c: A^{n-1}\to M$ such that $\partial c = \Ph\o d^n$.
This mapping $c$ induces a unique bimodule homomorphism 
$$\gather
\tilde \Ph: A\tilde\otimes(A/\Bbb R)^{\tilde\otimes(n-1)}\tilde\otimes A 
	\to M, \\
\tilde\Ph(a_0\otimes\bar a_1,\dots,\bar a_n, a_{n+1}) = 
	a_0\cdot c(a_1,\dots, a_n)\cdot a_{n+1}.
\endgather$$
and we have $\tilde \Ph \o I\o d^n = \tilde \Ph \o \partial \ph$, 
and moreover
$$\align
(\tilde \Ph \o \partial \ph)(a_1,\dots,a_n) 
&= \tilde\Ph(a_1\otimes\bar a_2\otimes\dots\otimes\bar a_n\otimes 1)\\
&\quad + \sum_{i=1}^{n-1}(-1)^i 
    	\tilde\Ph(1\otimes\bar a_1\otimes\dots\otimes\overline{a_ia_{i+1}}\otimes
	\dots\otimes\bar a_n\otimes 1)\\
&\quad +(-1)^n 
  	\tilde\Ph(1\otimes\bar a_1\otimes\dots\otimes\bar a_{n-1}\otimes a_n)\\
&= \partial c(a_1,\dots,a_n).
\endalign$$
So we get $\Ph\o d^n = \partial c = \tilde \Ph\o I\o d^n$ and the 
result follows from \nmb!{3.3}.

The second assertion of the proposition follows also from the last 
computation.
\qed\enddemo

\proclaim{\nmb.{3.6}. Corollary} For a convenient algebra $A$ and a 
convenient bimodule $M$ over $A$ we have 
$$H^n(A,M) \cong \frac{\Hom_A^A(\Om_n(A),M)}
{I^*(\Hom^A_A(A\tilde\otimes\bar 
A^{\tilde\otimes(n-1)}\tilde\otimes A,M))}.\qed$$
\endproclaim

\heading\totoc\nmb0{4}. The calculus of Fr\"olicher and 
Nijenhuis \endheading

\subheading{\nmb.{4.1}} In this section let $A$ be a convenient algebra 
with unit and let $\Om(A)=\Om_*(A)$ be the universal graded differential 
algebra for $A$.
The space $\operatorname{Der}_k\Om(A)$ consists of
all bounded \idx{\it (graded) derivations} of degree $k$, i.e.
all bounded linear mappings 
$D:\Om(A) \to  \Om(A)$ with $D(\Om_\ell(A)) \subset
\Om_{k+\ell}(A)$ and 
$D(\ph \ps) = D(\ph) \ps +(-1)^{k\ell}\ph
D(\ps)$ for $\ph \in \Om_\ell(A)$.
Obviously $\operatorname{Der}_k\Om(A)$ is a closed linear subspace
of $L(\Om (A),\Om (A))$ and thus a convenient vector space.

\proclaim{Lemma} The space 
$\Der\Om(A) = \bigoplus_k\Der_k\Om(A)$ is a
convenient graded Lie algebra with the graded commutator 
$[D_1,D_2] := D_1\o D_2 - (-1)^{k_1k_2}D_2\o D_1$ as bracket.
This means that the bracket is graded anticommutative,
$[D_1,D_2] = -(-1)^{k_1k_2}[D_2,D_1]$, and satisfies the graded
Jacobi identity
$$[D_1,[D_2,D_3]] = [[D_1,D_2],D_3] + (-1)^{k_1k_2}[D_2,[D_1,D_3]]$$
(so that $ad(D_1) = [D_1,\quad]$ is itself a derivation).
\endproclaim
\demo{Proof} Plug in the definition of the graded commutator and
compute. The boundedness of the bracket follows from \nmb!{1.11}.\qed\enddemo

\subheading{\nmb.{4.2}. Fields} 
Recall from \nmb!{2.6} that $d^*:\Hom^A_A(\Om_1(A),A)\to \Der(A;A)$ is an 
isomorphism, which we will also denote by $\L$.
We denote the space $\Hom^A_A(\Om_1(A),A)$ by $\X(A)$ and call it the 
space of \idx{\it fields} for the algebra $A$. Then 
$\L:\X(A)\to \Der(A;A)$ is an isomorphism of convenient vector spaces.
The space of derivations $\Der(A;A)$ is a convenient Lie algebra with the
commutator $[\quad,\quad]$ as bracket, and so we have an induced Lie bracket
on $\X(A)=\Hom^A_A(\Om_1(A),A)$ which is given by 
$\L([X,Y])a=[\L_X,\L_Y]a=\L_X\L_Ya-\L_Y\L_Xa$. It will be referred to 
as the Lie bracket of fields.

\proclaim{\nmb.{4.3}. Lemma} Each field 
$X\in \X(A)=\Hom^A_A(\Om_1(A),A)$ is by definition a bounded $A$-bimodule
homomorphism $\Om_1(A)\to A$. It prolongs uniquely to a graded derivation 
$j(X)=j_X:\Om(A)\to\Om(A)$ of degree $-1$ by
$$\align 
j_X(a)&=0\quad\text{ for }a\in A=\Om_0(A),\\
j_X(\om)&=X(\om)\quad\text{ for }\om\in\Om_1(A)\\
j_X(\om_1\otimes _A \dots\otimes _A \om_k)&=\\
    =\sum_{i=1}^{k-1}(-1)^{i-1}\om_1\otimes _A &\dots
    \otimes _A \om_{i-1}\otimes _A  X(\om_i)\om_{i+1}
    \otimes _A\dots\otimes _A \om_k\\
    +(-1)^{k-1}\om_1\otimes _A &\dots\otimes _A \om_{k-1}X(\om_k)
\endalign$$
for $\om_i\in\Om_1(A)$. The derivation $j_X$ is called the 
\idx{\it contraction operator} of the field $X$.
\endproclaim

\demo{Proof}
This is an easy computation
\qed\enddemo

With some abuse of notation we write also
$\om(X)=X(\om)=j_X(\om)$ for $\om\in\Om_1(A)$ and 
$X\in\X(A)=\Hom^A_A(\Om_1(A),A)$. 

\subheading{\nmb.{4.4}} A derivation $D \in \Der_k\Om(A)$ is called
\idx{\it algebraic} if $D\mid \Om_0(A) = 0$.
Then $D(a\om) =
aD(\om)$ and $D(\om a)=D(\om)a$ for $a\in A$, so $D$ restricts to a 
bounded bimodule homomorphism, an element of 
$\Hom^A_A(\Om_l(A),\Om_{l+k}(A))$.
Since we have
$\Om_l(A) = \Om_1(A)\tilde \otimes_A\dots\tilde\otimes_A\Om_1(A)$
and since for a product of one forms we have
$D(\om_1\otimes _A\dots\otimes _A\om_l)=\sum_{i=1}^l(-1)^{ik}
\om_1\otimes _A\dots\otimes _A D(\om_i)\otimes _A\dots\otimes_A\om_l$,
the derivation $D$ is uniquely determined by its restriction
$$K:=D|\Om_1(A)\in\Hom^A_A(\Om_1(A),\Om_{k+1}(A));$$
we write $D=j(K)=j_K$ to express this dependence. Note the defining 
equation $j_K(\om)= K(\om)$ for $\om\in\Om_1(A)$.
Since it will be very important in the sequel we will use the 
notation
$$\align\Om^1_k=\Om^1_k(A):&=\Hom^A_A(\Om_1(A),\Om_k(A))\\
\Om^1_*=\Om^1_*(A)&=\bigoplus_{k=0}^\infty\Om^1_k(A).\\
\endalign$$
Elements of the space $\Om^1_k$ will be called \idx{\it field valued 
$k$-forms}, those of $\Om^1_*$ will be called just \idx{\it field 
valued forms}.

\subheading{\nmb.{4.5}} In \nmb!{4.3} we have already 
met some algebraic graded derivations: for a
field $X\in\X(A)$ the derivation
$j_X$ is of degree $-1$.
The basic derivation $d$ is of
degree $1$. 
Note also that $\L_X := d\,j_X + j_X\,d$ translates to $\L_X = [j_X,d]$
and that this extends $\L_X$ from a derivation $A$ to a derivation of 
degree 0 of $\Om_*(A)$.

\proclaim{\nmb.{4.6} Theorem} 
\therosteritem1 
For $K \in \Om^1_{k+1}(A)$ and $\om_i\in\Om_1(A)$
the formula
$$j_K(\om_0\otimes _A  \dots\otimes _A  \om_\ell)
     =\sum_{i=0}^\ell(-1)^{ik}\om_0\otimes _A  \dots
     \otimes _A  K(\om_i)\otimes _A  \dots \otimes _A \om_k$$
defines an algebraic graded derivation $i_K \in \Der_k\Om(A)$ and any
algebraic derivation is of this form.

\therosteritem2
The map
$$j:\Om^1_{k+1}=\Hom^A_A(\Om_1(A),\Om_{k+1}(A))\to \Der^{alg}_k\Om(A)$$
where $\Der^{alg}_k\Om (A)$ denotes the closed linear subspace of
$\Der _k\Om (A)$ consisting of all algebraic derivations is an 
isomorphism of convenient vector spaces.

\therosteritem3
By $j([K,L]^{\De}):= [j_K,j_L]$ we get a bracket
$[\quad,\quad]^{\De}$ on the space 
$\Om^1_{*-1}$ which defines a convenient
graded Lie algebra structure with the grading as indicated, and
for $K\in \Om^1_{k+1}$,and 
$L\in \Om^1_{\ell+1}$ we have 
$$[K,L]^{\De} = j_K\o L - (-1)^{k\ell}j_L\o K.$$
\endproclaim
$[\quad,\quad]^{\De}$ is called the \idx{\it algebraic bracket} or
also the \idx{\it abstract De Wilde, Lecomte bracket} see
\cite{DeWilde, Lecomte, 1988}.

\demo{Proof} The first assertion is clear from the definition.

Clearly the map $D\mapsto D|\Om _1(A)$ is bounded. To show that $j$
is bounded recall that $\Der _d\Om (A)$ is a closed subspace of
$L(\Om (A),\Om (A))\cong \prod _kL(\Om _k(A),\Om (A))$. 
By \nmb!{2.7}.2 it suffices to show that $j$ is bounded as a map to
$L^A(\Om _1(A),\dots ,\Om _1(A);\Om (A))$ and by the linear 
uniform boundedness principle \nmb!{1.9}.2 it is enough to 
show that for all $\om _i\in \Om _1(A)$ the map $K\mapsto
j_K(\om _1\otimes _A  \dots \otimes _A  \om _k)$ is bounded. But this is 
clear by \therosteritem1.

For the third assertion it suffices to evaluate
$[j_K,j_L]$ at some $\om\in\Om_1(A)$.
\qed\enddemo

\subheading{\nmb.{4.7}} The exterior derivative $d$ is an element of
$\Der_1\Om(A)$. In view of the formula $\L_X = [j_X,d] = j_X\,d
+ d\,j_X$ for fields $X$, we define for $K \in
\Om^1_k$ the \idx{\it Lie derivation} 
$\L_K = \L(K) \in \Der_k\Om(A)$ by $\L_K := [j_K,d]$.

Then the mapping $\L:\Om^1_* \to  \Der\Om(A)$ is obviously bounded
and it is injective by the universal property of $\Om _1(A)$, since
$\L_Ka=j_Kda = K(da)$ for $a\in A$.

\proclaim{Theorem} For any graded derivation $D \in \Der_k\Om(A)$
there are unique homomorphisms 
$K \in \Om^1_k$ and 
$L \in \Om^1_{k+1}$
such that $$D = \L_K + j_L.$$
We have $L=0$ if and only if $[D,d]=0$. $D$ is algebraic if and
only if $K=0$.
\endproclaim

\demo{Proof} $D|A:a\mapsto Da$ is a derivation $A\to \Om_d(A)$, 
so by \nmb!{2.5} it is of the form $D|A=K\o d$ for a unique
$K \in \Om^1_k$. 

The defining equation for $K$ is  $Da =  j_Kda = \L_Ka$
for $a\in A$. Thus $D - \L_K$ is an algebraic
derivation, so $D - \L_K = j_L$ by \nmb!{4.4} for unique $L \in
\Om^1_{k+1}$. 

Since we have $[d,d] = 2d^2 =0$, by the graded Jacobi identity
we obtain $0 = [j_K,[d,d]] = [[j_K,d],d] + (-1)^{k-1}[d,[j_K,d]]
= 2[\L_K,d]$. The mapping $L \mapsto [j_L,d] = \L_L$ is
injective, so the last assertion follows.
\qed\enddemo

\subheading{\nmb.{4.8}. The Fr\"olicher-Nijenhuis bracket} 
Note that $j(Id_{\Om_1(A)})\om =
k\om$ for $\om \in \Om_k(A)$. Therefore we have 
$\L(Id_{\Om_1(A)})\om = j(Id_{\Om_1(A)})d\om - d\,j(Id_{\Om_1(A)})\om =
(k+1)d\om - kd\om = d\om$. Thus $\L(Id_{\Om_1(A)}) = d$.

\subheading{\nmb.{4.9}} 
Let $K \in \Om^1_k$ and 
$L \in \Om^1_\ell$. 
Then obviously $[[\L_K,\L_L],d] =0$, so we have
$$[\L(K),\L(L)] = \L([K,L])$$
for a uniquely defined 
$[K,L] \in \Om^1_{k+\ell}$. This
vector valued form $[K,L]$ is called the \idx{\it abstract
Fr\"olicher-Nijenhuis bracket} of $K$ and $L$.

\proclaim{Theorem} The space
$\Om^1_* = \bigoplus_k\Om^1_k$
with its usual grading and the Fr\"olicher-Nijen\-huis 
bracket is a convenient graded Lie algebra.
$Id_{\Om_1(A)} \in \Om^1_1$ is in the center, 
i.e. $[K,Id_{\Om_1(A)}] = 0$ for all $K$.

$\L:(\Om^1_*, [\quad,\quad]) \to  \Der\Om(A)$ is a bounded injective
homomorphism of graded Lie algebras. For fields in 
$\Hom^A_A(\Om_1(A),A)$, i\. e\. bounded derivations of $A$,
the Fr\"olicher-Nijenhuis bracket coincides with the bracket defined in 
\nmb!{4.2}.
\endproclaim
\demo{Proof}
Boundedness of the bracket follows from the fact that the map 
$\L_K\mapsto K$ is bounded as it is just the composition of the 
restriction to $A$ with the bounded inverse to $d^*$ constructed in
\nmb!{2.6}.

For $X,Y\in \Hom^A_A(\Om_1(A),A)$ we have
$j([X,Y])da = \L([X,Y])a = [\L_X,\L_Y]a$. The rest
is clear. \qed\enddemo

\proclaim{\nmb.{4.10}. Lemma} 
For homomorphisms $K \in \Om^1_k$ and 
$L \in \Om^1_{\ell+1}$ we have
$$\align [\L_K,j_L] &= j([K,L]) - (-1)^{k\ell}\L(j_L\o K)\text{, or}\\
[j_L,\L_K] &= \L(j_L\o K) - (-1)^k\,j([L,K]).\endalign$$
\endproclaim

\demo{Proof} For $a \in A$ we have 
$[j_L,\L_K]a = j_L\,j_K\,da - 0 = j_L(K(da)) = (j_L\o K)(da) =
\L(j_L\o K)a$. 
So $[j_L,\L_K] - \L(j_L\o K)$ is an algebraic derivation.
$$\multline 
[[j_L,\L_K],d] = [j_L,[\L_K,d]] -(-1)^{k\ell}[\L_K,[j_L,d]] = \\
= 0 -(-1)^{k\ell}\L([K,L]) = (-1)^k[j([L,K]),d]).
\endmultline$$
Since $[\quad,d]$ kills the `$\L$'s' and is injective on the
`$j$'s', the algebraic part of $[j_L,\L_K]$ is $(-1)^k\,j([L,K])$.
\qed\enddemo

\proclaim{\nmb.{4.11}. Theorem} For homomorphisms
$K_i \in \Om^1_{k_i}$ and 
$L_i \in \Om^1_{k_i+1}$ we have
$$\align
&\left[\L_{K_1}+j_{L_1},\L_{K_2}+j_{L_2}\right] = \tag1 \\
&\qquad = \L\left([K_1,K_2] 
    + j_{L_1}\o K_2 - (-1)^{k_1k_2}j_{L_2}\o K_1\right)\\
&\qquad\quad + i\left([L_1,L_2]^{\De} 
    + [K_1,L_2] -(-1)^{k_1k_2}[K_2,L_1]\right).
\endalign$$
Each summand of this formula looks like a semidirect product of
graded Lie algebras, but the mappings
$$\align 
j:\Om^1_{*-1} &\to \End(\Om^1_*,[\quad,\quad]) \\
\ad:\Om^1_* &\to \End(\Om^1_{*-1},[\quad,\quad]^{\De}),
     \quad \ad_KL=[K,L],
\endalign$$
do not take values in the subspaces of graded derivations. We
have instead for homomorphisms 
$K \in \Om^1_k$ and 
$L \in \Om^1_{\ell+1}$ 
the following relations:
$$\align &j_L\o [K_1,K_2] = [j_L\o K_1,K_2] 
     + (-1)^{k_1\ell}[K_1,j_L\o K_2] \tag2 \\
&\qquad -\left((-1)^{k_1\ell}j(\ad_{K_1}L)\o K_2 
      - (-1)^{(k_1+\ell)k_2}j(\ad_{K_2}L)\o K_1\right)\\ \allowdisplaybreak
&\ad_K[L_1,L_2]^{\De} = [\ad_KL_1,L_2]^{\De} 
      + (-1)^{kk_1}[L_1,\ad_KL_2]^{\De} -\tag3 \\
&\qquad -\left((-1)^{kk_1}\ad(j(L_1)\o K)L_2 
      - (-1)^{(k+k_1)k_2}\ad(j(L_2)\o K)L_1\right)
\endalign$$
\endproclaim

The algebraic meaning of the relations of this theorem and its
consequences in group theory have been investigated in
\cite{Michor, 1990}. The corresponding product of groups is well
known to algebraists under the name `Zappa-Szep'-product.

\demo{Proof} Equation \thetag1 is an immediate consequence of
\nmb!{4.10}. Equations \thetag2 and \thetag3 follow from \thetag1
by writing out the graded Jacobi identity, or as follows:
Consider $\L(j_L\o [K_1,K_2])$ and use \nmb!{4.10} repeatedly to
obtain $\L$ of the right hand side of \thetag2. Then consider
$j([K,[L_1,L_2]^{\De}])$ and use again \nmb!{4.10} several times to
obtain $i$ of the right hand side of \thetag3.
\qed\enddemo

\subheading{\nmb.{4.12}. Naturality of the Fr\"olicher-Nijenhuis
bracket} Let $f:A \to  B$ be a bounded algebra homomorphism.
Two forms 
$K \in \Om^1_k(A)=\Hom^A_A(\Om_1(A),\Om_k(A))$ and 
$K' \in \Om^1_k(B) = \Hom^B_B(\Om_1(B),\Om_k(B))$ 
are called \idx{\it $f$-related} or \idx{\it
$f$-dependent}, if we have
$$K'\o \Om_1(f) = \Om_k(f)\o K: \Om_1(A)\to\Om_k(B) ,\tag{1}$$
where $\Om_*(f)$ is described in \nmb!{2.10}.

\proclaim{Theorem}
\roster
\item[2] If $K$ and $K'$ as above are $f$-related then
  $j_{K'}\o \Om(f) = \Om(f) \o j_{K}:\Om(A) \to  \Om(B)$.
\item If $j_{K'}\o \Om(f)| d(A) = \Om(f)\o j_{K}| d(A)$, then 
  $K$ and $K'$ are $f$-related, where $d(A)\subset \Om_1(A)$ 
  denotes the space of exact $1$-forms.
\item If $K_j$ and $K_j'$ are $f$-related for $j=1,2$, then 
  $j_{K_1}\o K_2$ and $j_{K'_1}\o K'_2$ are $f$-related, and also
  $[K_1,K_2]^{\De}$ and $[K'_1,K'_2]^{\De}$ are $f$-related.
\item If $K$ and $K'$ are $f$-related then 
  $\L_{K'}\o \Om(f) = \Om(f)\o \L_{K}:\Om(A) \to  \Om(B)$.
\item If $\L_{K'}\o \Om(f)\mid \Om_0(A) = \Om(f)\o\L_{K}\mid \Om_0(A)$, 
  then $K$ and $K'$ are $f$-related.
\item If $K_j$ and $K_j'$ are $f$-related for $j=1,2$, then
  their Fr\"olicher-Nijenhuis brackets $[K_1,K_2]$ and
  $[K'_1,K'_2]$ are also $f$-related.
\endroster\endproclaim

\demo{Proof} \therosteritem2. Since both sides are graded derivations 
over $\Om(f)$ it suffices to check this for a 1-form 
$\om\in \Om_1(A)$.
By \nmb!{4.6} and \nmb!{2.10} we have
$\Om_{k}(f)j_K(\om)  = \Om_k(f)K(\om) = K'(\Om_1(f)\om) = 
     j_{K'}\Om_1(f)(\om)$.

\therosteritem3 follows from the universal 
property of $\Om _1(A)$ because $K'\o \Om _1(f)\o d$ and
$\Om _k(f)\o K\o d$ are both derivations from $A$ into $\Om _k(B)$
which is an $A$-bimodule via $f$ and the multiplication in 
$\Om (B)$.

\therosteritem4 is obvious;
the result for the bracket then follows from \nmb!{4.6}.3.

\therosteritem5 The algebra homomorphism $\Om(f)$ intertwines the
operators $j_K$ and $j_{K'}$ by \therosteritem2, and $\Om(f)$
commutes with the exterior derivative $d$. Thus $\Om(f)$
intertwines the commutators $[j_K,d] = \L_K$ and $[j_{K'},d] =
\L_{K'}$.

\therosteritem6 For an element $g \in \Om_0(A)$ we have 
$\L_K\,\Om(f)\,g = j_K\,d\,\Om(f)\,g = j_K\,\Om(f)\,dg$ and 
$\Om(f)\,\L_{K'}\,g = \Om(f)\,j_{K'}\,dg$. By \therosteritem3 the
result follows.

\therosteritem7 The algebra homomorphism $\Om(f)$ intertwines 
$\L_{K_j}$ and $\L_{K'_j}$, so also their graded commutators
which equal  $\L([K_1,K_2])$ and $\L([K'_1,K'_2])$, respectively.
Now use \therosteritem6 . \qed\enddemo

\heading\totoc\nmb0{5}. Distributions and integrability\endheading

\subheading{\nmb.{5.1}. Distributions} By a \idx{\it distribution} 
in a convenient algebra $A$ we mean a $c^\infty$-closed sub-$A$-bimodule
$\Cal D$ of $\Om_1(A)$.

The distribution $\Cal D$ is called \idx{\it globally integrable} if 
there exists a $c^\infty$-closed subalgebra $B$ of $A$ such that
$\Cal D$ is the $c^\infty$-closure in $\Om _1(A)$ of the subspace
generated by $A(d(B))$ and $d(B)A$.

The distribution $\Cal D$ is called \idx{\it splitting} if there 
exists a bounded projection $P\in\Om^1_1(A)=\Hom^A_A(\Om_1(A),\Om_1(A))$
onto $\Cal D$, i.e. $P\o P=P$ and $\Cal D = P(\Om_1(A))$.  Then there 
is a complementary submodule $\ker P\subset \Om_1(A)$.

The distribution $\Cal D$ is called \idx{\it involutive} if the 
$c^\infty$-closed ideal $(\Cal D)_{\Om_*(A)}$ generated by $\Cal D$ in the
graded algebra $\Om_*(A)$ is stable under $d$, i\.e\. if
$d(\Cal D)\subset (\Cal D)_{\Om_*(A)}$.

\subheading{\nmb.{5.2}. Comments}
One should think of this as follows: In differential geometry, where 
we have $A = C^\infty(M,\Bbb R)$ for a manifold $M$,
a distribution is usually given as a sub vector bundle $E$ of the 
tangent bundle $TM$. Then
$\Cal D$ is the $A$-bimodule of those 1-forms which annihilate 
the subbundle $E$ of $TM$. Global integrability then means that it is 
integrable and that the space of functions which are constant along 
the leaves of the foliation generates those forms. This is a strong 
condition: There are foliations where this space of functions 
consists only of the constants, and this can be embedded into any 
manifold. So in $C^\infty(M,\Bbb R)$ there are always involutive 
distributions which are not globally integrable. To prove some
Frobenius theorem a notion of local integrability would be necessary.

\subheading{\nmb.{5.3} Curvature and cocurvature}  
Let $P\in\Om^1_1(A)=\Hom^A_A(\Om_1(A),\Om_1(A))$
be a projection, then the image $P(\Om_1(A))$ is a splitting distribution, 
called the \idx{\it vertical distribution} of $P$ and 
the complement $\ker P$ is also a splitting distribution, called the 
\idx{\it horizontal\ign{ distribution}} one.
$\bar P:= Id_{\Om_1(A)}- P$ is a projection onto  the 
horizontal distribution.

We consider now the Fr\"olicher-Nijenhuis bracket $[P,P]$ of $P$ and 
define 
$$\alignat2
R &= R_P = [P,P]\o P &\qquad &\text{ the \it curvature,}\\
\bar R &= \bar R_P = [P,P]\o \bar P &\qquad &\text{ the \it cocurvature.}
\endalignat$$
The \idx{curvature} and the \idx{cocurvature} are elements of   
$\Om_2^1(A)=\Hom^A_A(\Om_1(A),\Om_2(A))$. The curvature kills 
elements of the horizontal distribution, so it is \idx{\it vertical}.
The cocurvature kills elements of the vertical distribution.

Since the identity $Id\in\Om^1_1(A)$ lies in the center of the 
Fr\"olicher-Nijenhuis algebra we get 
$[\bar P,\bar P]=[Id-P,Id-P] = [P,P]$ and hence 
$\bar R_P=R_{\bar P}$.
We shall also need the homomorphisms of graded algebras 
$\Om(P), \Om(\bar P):\Om(A) \to \Om(A)$ with 
$\Om_0(P)=\Om_0(\bar P) = Id_A$ which are induced by the 
bimodule homomorphisms $P, \bar P: \Om_1(A)\to \Om_1(A)$.

\proclaim{\nmb.{5.4}. Lemma} In the setting of \nmb!{5.3} the 
following assertions hold:

1. For $\om\in\Om_1(A)$ we have 
$$\align
R_P(\om) &= [P,P](P(\om)) = -2(\Om(\bar P)\o d\o P)(\om) \\
\bar R_P(\om) &= [P,P](\bar P(\om)) = -2(\Om(P)\o d\o \bar P)(\om). 
\endalign$$

2. For the $c^\infty$-closed ideals generated by the distributions
$\ker P$ and $P(\Om_1(A))$ we have  $(\ker P)_{\Om_*(A)}= \ker \Om(P)$ and 
$(P(\Om_1(A)))_{\Om_*(A)} = \ker \Om(\bar P)$.

3. The curvature $R=[P,P]\o P$ is zero if and only if the horizontal 
distribution is involutive. The cocurvature $\bar R=[P,P]\o (Id-P)$ 
is zero if and only if the vertical distribution
$P(\Om_1(A))$ is involutive.  
\endproclaim

\demo{Proof}
\therosteritem1 
It suffices to show the first equation. For $\om\in \Om_1(A)$ we 
have:
$$\align
[P,P](\om)&= [\bar P,\bar P](\om) = j([\bar P,\bar P])(\om) \\
& = [\L_{\bar P},j_{\bar P}](\om) + \L(j_{\bar P}\bar P)(\om)  
     \quad\text{ by \nmb!{4.10}}  \\
& = \L_{\bar P}j_{\bar P}(\om) - j_{\bar P}\L_{\bar P}(\om) + \L_{\bar P}(\om)
     \quad\text{ since }j_{\bar P}{\bar P}={\bar P}^2={\bar P}\\
& = 2(j_{\bar P}d{\bar P}(\om)-d{\bar P}(\om)) - 
     j_{\bar P}j_{\bar P}d(\om) + j_{\bar P}d(\om).
\endalign$$
For $\om$, $\ph\in\Om_1(A)$ we have 
$$\align
j_{\bar P}j_{\bar P}(\om\otimes _A  \ph) &= 
     j_{\bar P}({\bar P}(\om)\otimes _A \ph + \om\otimes _A
     {\bar P}(\ph))\\
&= {\bar P}(\om)\otimes _A  \ph + 2 {\bar P}(\om)\otimes _A
     {\bar P}(\ph) 
     + \om\otimes _A {\bar P}(\ph)\\
&= (2\Om({\bar P}) + j_{\bar P})(\om\otimes _A \ph),\quad\text{ thus }\\
j_{\bar P}j_{\bar P}|\Om_2(A) &= (2\Om({\bar P})+j_{\bar P})|\Om_2(A).
\endalign$$
So we have 
$$\align
[P,P](\om)&= 2(j_{\bar P}d{\bar P}(\om)
     -d{\bar P}(\om)-\Om({\bar P})(d(\om)))\\
R_P(\om) &= [P,P](P(\om))= -2\Om({\bar P})dP(\om)
\endalign$$
as required.

\therosteritem2 The kernel of the bounded algebra homomorphism $\Om(P)$ is a
$c^\infty$-closed ideal and contains $\ker P$. On the other hand any 
$\om\in\Om_1(A)\otimes _A\dots \otimes _A\Om _1(A)\cap \ker\Om(P)$
(non-completed tensor product) may be written as a finite sum 
$\om=\sum_i\om_{1,i}\otimes _A\dots\otimes _A\om_{k,i}$ with the property 
that 
$\sum_iP(\om_{1,i})\otimes _A\dots\otimes _AP(\om_{k,i})=0$. Since 
$P+\bar P=Id_{\Om_1(A)}$ we have $\om_{j,i} = 
P(\om_{j,i})+\bar P(\om_{j,i})$ for all $j$ and $i$. Thus each 
summand of $\om$ splits into a sum of products of  $P(\om_{j,i})$  and 
$\bar P(\om_{j,i})$ and the sum of those products containing only 
$P(\om_{j,i})$ vanishes. So at least one $\bar P(\om_{j,i})$ appears 
in each summand and the whole sum is in the ideal generated by 
$\ker \Om _1(P) = \bar P(\Om_1(A))$.

By \nmb!{1.7} $\Om _k(A)\cap \ker (\Om (P))$ is the completion
of $\Om_1(A)\otimes _A\dots \otimes _A\Om _1(A)\cap \ker\Om(P)$
so it must be the  $c^\infty$-closure in $\Om _k(A)$ of this space
and hence must also be contained in the $c^\infty$-closed ideal.

The second assertion follows by symmetry.

\therosteritem3 We have to prove only the first assertion. The 
distribution $\ker\bar P$ is 
involutive if and only if for all $\om\in\Om_1(A)$ we have 
$dP\om \in (\ker P)_{\Om_*(A)} = \ker \Om(P)$. By \therosteritem2 this 
is equivalent to $R(\om)= -2\Om(\bar P)(dP(\om))=0$ for all 
$\om\in\Om_1(A)$.
\qed\enddemo

\proclaim{\nmb.{5.5}. Lemma (Bianchi identity)} If 
$P \in \Om^1_1(A)$ is a
projection with curvature $R$ and cocurvature
$\bar R$, then we have
$$\align &[P ,R+\bar R] = 0 \\
         &2[R,P ] = j_R\bar R + j_{\bar R}R.
\endalign $$

\endproclaim
\demo{Proof}  We have $[P ,P ] = R +\bar R$ by \nmb!{5.3} and 
$[P ,[P ,P ]] = 0$ by the graded Jacobi identity. So the first
formula follows. We have $R =[P ,P ]\o P = j_{[P ,P ]}\o P $.
By \nmb!{4.11}.2 we get 
$j_{[P ,P ]}\o [P ,P ]=2[j_{[P ,P ]}\o P ,P ]-0=2[R,P ]$. 
Therefore 
$2[R,P ] = j_{[P ,P ]}\o [P ,P ] = j(R+\bar R)\o (R+\bar R) =
j_R\o\bar R + j_{\bar R}\o R$ 
since $R$ has vertical values and kills
vertical vectors, so $j_R\o R = 0$; likewise for $\bar R$.
\qed\enddemo

\heading\totoc\nmb0{6}. Bundles and connections \endheading

Let $G$ be a Lie group in the usual sense. 
We want to carry over to non-commutative /n differential geometry the concepts of principal bundles, 
characteristic classes, and Chern-Weil homomorphism. The last two 
concepts still make difficulties, since we do not know how to express 
local triviality and only some of the usual properties hold in the 
general setup we use. 

\subheading{\nmb.{6.1}. Definition} 
By a bundle in
non-commutative differential geometry we mean
a convenient algebra $A$ together with a closed subalgebra 
$B\hookrightarrow A$.

The bundle is said to have a finite dimensional Lie group $G$ 
as \idx{\it structure group} if we have an injective
homomorphism $\la: G\to \operatorname{Aut}(A)$, such that $\la:G \to
L(A,A)$ is smooth
and $B=A^G$, the subalgebra of all
elements fixed by the $G$-action. 

We remark that for the notion of a principal bundle one should 
add requirements like quantum transitiveness on the fiber, compare 
with \cite{Narnhofer, Thirring, Wicklicky, 1988}, but this is still not 
enough to get the Chern-Weil homomorphism, see also \nmb!{6.9}.

If $p:P\to M$ is a smooth principal bundle in the usual sense, we put
$A=C^\infty(P,\Bbb R)$ and $B=C^\infty(M,\Bbb R)$, which is embedded
into $A$ via $p^*$. Then clearly all requirements are satisfied.

\proclaim{\nmb.{6.2}. Lemma} For each $g\in G$ the algebra
automorphism $\la_g:A\to A$ extends to an automorphism of the algebra
of differential forms as follows:
$$\CD
A @>>> \Om(A)\\
@V\la_gVV  @V\la_gVV\\
A @>>> \Om(A).
\endCD$$
\endproclaim

\demo{Proof}
This follows from the universal property \nmb!{2.9}.
\qed\enddemo

\subheading{\nmb.{6.3}. Horizontal forms}
Recall, that on a classical bundle the horizontal forms are
exactly those which annihilate vertical vectors. Guided by this we
define the space of \idx{\it horizontal 1-forms} $\Om^{\text{hor}}_1(A)$
as the closed $A$-bimodule generated by $\Om_1(B)$ in $\Om_1(A)$, in 
the bornological topology. Likewise
we define the algebra $\Om^{\text{hor}}(A)$ of all \idx{\it horizontal
forms} as the 
closed subalgebra of $\Om(A)$ generated by
$A+\Om(B)$.

So $\Om^{\text{hor}}_1(A)$ is the closed linear subspace generated by 
all elements of the form 
$a(db)a'$ for $a$, $a'\in A$ and $b\in B$. Since in 
$\Om_1(A)\subset A\tilde\otimes A$ we have 
$a(db)a'=a(1\otimes b-b\otimes1)a'=a\otimes ba'-ab\otimes a'$, we get
$A\tilde\otimes A/\Om^{\text{hor}}_1(A)=A\tilde\otimes_B A$ where $A$ 
is viewed as a $B$-bimodule. The situation is explained in the  
following diagram  
$$\CD
	@.	0				  @.		0                     @.	  0	@.\\
@.		@VVV						@VVV					  @VVV	@.\\
0	@>>>	\Om^{\text{hor}}_1(A) @=		\Om^{\text{hor}}_1(A) @>>> 0 @.\\
@.		@VVV						@VVV					  @VVV	@.\\	
0	@>>> \Om_1(A)			  @>>>	A\tilde\otimes A	  @>\mu>> A	@>>> 0\\
@.		@VVV						@VVV					  @VVV	@.\\	
0	@>>> \Om_1(A)/\Om^{\text{hor}}_1(A)  @>>> A\tilde\otimes_B A  @>\mu>> A	@>>> 0\\
@.		@VVV						@VVV					  @VVV	@.\\	
	@.	0				  @.		0                     @.	  0	@.
\endCD$$ 
which has exact columns and also rows since the middle row is splitting.

\subheading{\nmb.{6.4}. Principal connections} We have a good
description of horizontal forms, whereas vertical vector fields do
not exist in sufficient supply, thus we describe connections in the
form of horizontal projections. So a \idx{\it connection} on a bundle
$B\hookrightarrow A$ is an element 
$\ch\in\Om^1_1(A) = \Hom^A_A(\Om_1(A),\Om_1(A))$ which satisfies
$\ch\o\ch=\ch$ (equivalently $j_\ch\o\ch=\ch$), such that the image
of $\ch$ is $\Om^{\text{hor}}_1(A)$, the space of horizontal 1-forms
of the bundle.

Note that a connection $\ch:\Om_1(A)\to \Om^{\text{hor}}_1(A)$ has a
unique extension as an $A$-bimodule homomorphism 
$$\gather
\Om_k(A)= \Om_1(A)\otimes_A\dots\otimes_A\Om_1(A) @>\Om(\ch)>>
     \Om^{\text{hor}}_1(A)\otimes_A\dots\otimes_A\Om^{\text{hor}}_1(A)\\
\om_1\otimes \dots \otimes\om_k\mapsto 
     \ch(\om_1)\otimes \dots \otimes\ch(\om_k).
\endgather$$

A connection $\ch$ on a bundle with structure group $G$ 
is called a \idx{\it principal connection} if it is $G$-equivariant: 
$\ch\o\la_g=\la_g\o\ch$ for all $g\in G$.

For a usual principal bundle this corresponds to the projection of 
forms onto horizontal forms, which describe the vertical 
distribution. This explains our choice of names here and in 
\nmb!{5.3}.

PROBLEM: What means `locally trivial' for a bundle? Does it imply the
existence of connections?
    
\subheading{\nmb.{6.5}. Curvature} Let $\ch$ be a connection on 
a non-commutative bundle $B\hookrightarrow A$. The \idx{\it curvature}
$R=R(\ch)$ of the connection is given by
$$R=[\ch,\ch]\in \Om^1_2(A)=\Hom^A_A(\Om_1(A),\Om_2(A)),$$
the abstract Fr\"olicher-Nijenhuis bracket of $\ch$ with itself.

\proclaim{\nmb.{6.6}. Lemma} The curvature of a connection satisfies
$$R\in\Hom^A_A\left(\Om_1(A)/
	\Om_1^{\text{hor}}(A),\Om_2^{\text{hor}}(A)\right).$$
If the connection is principal then also $R$ is $G$-equivariant.
\endproclaim

\demo{Proof}
By definition $\Om_1^{\text{hor}}(A) = \ch(\Om_1(A))$ is globally 
integrable, thus $R_\ch = [\ch,\ch]\o\ch = 0$ and we have
$$\alignat2
R :&= [\ch,\ch] = \bar R_\ch = [\ch,\ch]\o (Id-\ch)&\quad&\text{ by 
     \nmb!{5.3}} \\
&= -2\Om(\ch)\o d \o (Id-\ch)&\quad&\text{ by \nmb!{5.4}.1}.
\endalignat$$
The last expression implies the first assertion. If $\ch$ is a 
principal connection it is $G$-equivariant and by \nmb!{4.12} also
$R=[\ch,\ch]$ is $G$-equivariant.
\qed\enddemo

\subheading{\nmb.{6.7}. Steps towards the Chern-Weil homomorphism}
Let $B\subset A$ be a non-commutative bundle with structure group 
$G$. Let $\frak g$ denote the Lie algebra of $G$. 
We differentiate the action $\la: G\to \Aut(A)$ and get bounded 
linear mappings
$$\CD
\frak g @>T_e\la>> \Der(A;A)\\
@|                 @AAd^*A    \\
\frak g @>\la'>>  \Hom^A_A(\Om_1(A),A).
\endCD$$
Using this we define a mapping 
$$\gather
\al: \Om_1(A)/\Om_1^{\text{hor}}(A) \to A\otimes \frak g^* \\
(Id_A\otimes \ev_X)\al(\om):= \la'(X)(\om)\text{ for }X\in\frak g, 
     \om\in\Om_1(A).
\endgather$$

\proclaim{\nmb.{6.8}. Lemma} This mapping $\al$ is well defined, 
an $A$-bimodule homomorphism, and is 
$G$-equivariant for the action $\la_g\otimes \Ad(g\i)^*$ on the right 
hand side.
\endproclaim

\demo{Proof}
For $X\in\frak g$, $a, a'\in A$, and $\om\in\Om_1(A)$  we have 
$$\align
(A\otimes \ev_X)\al(a\om a') &= \la'(X)(a\om a')  \\
&= a\la'(X)(\om)a'\quad\text{ since }\la'(X)\in \Hom^A_A(\Om_1(A),A)\\
&= (A\otimes \ev_X)(a\al(\om)a'),
\endalign$$
so $\al$ is a bimodule homomorphism. For $b\in B$ we have
$$\align
(A\otimes \ev_X)\al(a(db) a') &= a\la'(X)(db)a'\\
&= a(T_e\la.X)(b)a' = 0\quad\text{ since }\la_g(b)=b.
\endalign$$
So $\al$ annihilates horizontal forms and is thus well defined.
In order to prove that $\al$ is $G$-equivariant we begin with the 
following computation, where $g\in G$:
$$\align
\la_g(T_e\la.X)(a) &= \la_g \tfrac d{dt}|_0\la_{\exp tX}(a)\\
&= \tfrac d{dt}|_0 \la_g\la_{\exp tX}(a)
     \quad\text{ since }\la_g\text{ is linear and bounded}\\
&= \tfrac d{dt}|_0 \la_{g\exp(tX)\,g\i}(\la_g(a))\\
&= T_e\la(\Ad(g)X)(\la_g(a)).
\endalign$$
By the universality of $d$ we have $\Om_1(\la_g)\o d = d\o \la_g$ and 
thus we get
$$\align
\la_g(\la'(X)(a\,da')) &= \la_g(a\la'(X)(da')) = \la_g(a)\la_g(\la'(X)da')\\
&= \la_g(a)\,\la_g((T_e\la.X)(a')) \\
&= \la_g(a)\,(T_e\la.\Ad(g)X)(\la_g(a')) \\
&= \la_g(a)\,\la'(\Ad(g)X))(d\la_g(a')) \\
&= \la'(\Ad(g)X))(\la_g(a\,da')).
\endalign$$
So finally we have
$$\align
(A\otimes \ev_X)\al(\Om_1(\la_g)\om) &= \la'(X)(\Om_1(\la_g)\om)\\
&= \la_g(\la'(\Ad(g\i)X)\om)\\
&= (\la_g\otimes \ev_{\Ad(g\i)X})\al(\om)\\
&= (A\otimes \ev_X)(\la_g\otimes \Ad(g\i))\al(\om),
\endalign$$
so $\al\o\Om_1(\la_g) = (\la_g\otimes \Ad(g\i))\o \al$ as required.
\qed\enddemo

\subheading{\nmb.{6.9}. Remarks} We stop our development here and add 
just some remarks about the Chern-Weil homomorphism. To continue from 
this point one should add requirements to the bundle $A$ which imply 
that $\al$ is invertible (the inverse then describes the fundamental 
vector field mapping) and that the extension of the inverse to 
invariant polynomials on $\frak g$ factors to the $\bar\Om(B)$. 

A good model for the Chern-Weil homomorphism
is described in the paper \cite{Lecomte, 1985} where the following 
construction is given: 

Let $P\to M$ be a smooth principal fiber bundle with structure group 
$G$. Then we have the following exact sequence of vector bundles over 
$M$:
$$ 0 \to P[\frak g, Ad]\to TP/G @>Tp>> TM \to 0.$$
The smooth sections of these bundles give rise to the following exact 
sequence of Lie algebras:
$$ 0\to \X_{\text{vert}}(P)^G \to \X_{\text{proj}}(P)^G 
	\to \X(M)\to 0,$$
namely first all vertical $G$-equivariant vector fields (the Lie 
algebra of the gauge group), second the all projectable 
$G$-equivariant vector fields on $P$ (the infinitesimal principal 
bundle automorphisms), third all vector fields on the base.
The `dual' of this sequence of Lie algebras is
$$0\gets (\Om_*(A)/\Om^{\text{hor}}(A))^G \gets \Om_*(A)^G 
\gets \Om_*(B)\gets 0,$$
where $A=C^\infty(P,\Bbb R)$ and $B=C^\infty(M,\Bbb R)$.
For general algebras this sequence is not exact.
For any short exact sequence of Lie algebras \cite{Lecomte, 1985} has 
described a generalization of the Chern-Weil homomorphism in purely 
algebraic terms, using Chevalley cohomology of the Lie algebras in 
question. This should be the starting point of the Chern-Weil 
homomorphism in non-commutative differential geometry.    

\heading\totoc\nmb0{7}. Polyderivations and the Schouten-Nijenhuis 
bracket\endheading

In this section we describe the analogue of the Schouten-Nijenhuis 
bracket in the setting of non-commutative differential geometry. It 
turns out that one has to require skew symmetry in the construction 
in order to get a meaningful theory. In the end we obtain the Poisson 
structures for convenient algebras. The results in this section are 
also a generalization for non-commutative algebras of the results in 
\cite{Krasil'shchik, 1988}, which were the original motivation for 
the developments here, but our approach is different:
we first show that the Nijenhuis-Richardson bracket
(c\.f\. \cite{Nijenhuis, Richardson, 1967} and
\cite{Lecomte, Michor, Schicke\-tanz})
passes to the convenient setting and then by restricting it to a 
suitable space of polyderivations (the non-commutative analog of
multi vector fields) we derive a generalization of the Schouten-Nijenhuis
bracket.

\subheading{\nmb.{7.1}} It has been noticed in 
\cite{De Wilde, Lecomte, 1985} that for any  
smooth manifold $M$ the Schouten-Nijenhuis bracket on the space 
$C^\infty(\La TM)$ of 
all multivector fields imbeds as a graded sub Lie algebra into the 
space $\La^{*+1}(C^\infty(M,\Bbb R);C^\infty(M,\Bbb R))$ with the 
Nijenhuis-Richardson bracket (see \nmb!{7.2} for a description of 
this space). 
Lecomte told us, that a very elegant 
proof of this fact can be given in the following way: The space 
$C^\infty(M,\Bbb R)$ of smooth functions is the degree $-1$  part of 
the Schouten-Nijenhuis algebra. By the universal property of the 
Nijenhuis-Richardson algebra 
$(\La^{*+1}(C^\infty(M,\Bbb R);C^\infty(M,\Bbb R)),[\quad,\quad]^\wedge )$ 
described in 
\cite{Lecomte, Michor, Schicketanz} the identity on 
$C^\infty(M,\Bbb R)$ prolongs to a unique homomorphism $\Ph$ of graded Lie 
algebras from the Schouten-Nijenhuis algebra into the 
Nijenhuis-Richardson algebra, and a simple computation described in 
\cite{Lecomte, Melotte, Roger, 1989} shows that 
$\Ph(T)=d^*(T)=T\o(d\x \dots\x d)$, where 
$d$ is the exterior differential.

This shows that the Schouten-Nijenhuis bracket which we will 
construct below boils down to the usual one in the commutative case 
$A=C^\infty(M,\Bbb R)$.

\subheading{\nmb.{7.2} The Nijenhuis-Richardson bracket in the 
convenient setting} Let $V$ be a convenient vector space. 
We consider the space $\La^k(V)$ of all 
bounded $k$-linear skew symmetric functionals 
$V\x\dots\x V\to \Bbb R$, where 
$\La^0(V)=\Bbb R$. Then $\La(V)=\bigoplus_{k\ge 0} \La^k(V)$ is a 
graded commutative convenient algebra with the usual wedge product
$$\multline
(\ph\wedge\ps)(\row v1{k+\ell}) = \\
= \tfrac1{k!\ell!}\sum_\si \sign\si\; 
\ph(\row v{\si1}{\si k})\ps(\row v{\si(k+1)}{\si(k+\ell)}),
\endmultline \tag 1$$
where the sum is over all permutation of $k+\ell$ symbols.

Now let $W$ be another convenient vector space.
We need the space $\La^k(V;W)$ of all bounded $k$-linear 
mappings $V\x\dots\x V\to W$. Then 
$\La(V;W)=\bigoplus_{k\ge0}\La^k(V,W)$ is a graded convenient vector 
space and a graded convenient module over the graded commutative 
algebra $\La(V)$ with the wedge product \thetag1 from above.   
If $A$ is a convenient algebra then $\La(V;A)$ is an associative 
graded convenient algebra with the (formally) same wedge product.

Now for $K\in\La^{k+1}(V;V)$ and $\Ph\in\La^p(V;W)$ we define 
$$\multline
(i_K\Ph)(\row v1{k+p}) = \\
= \tfrac1{(k+1)!(p-1)!}\sum_\si \sign\si\;
\Ph(K(\row v{\si1}{\si(k+1)}),\row v{\si(k+2)}{\si(k+p)}).
\endmultline\tag2$$
Then the following results hold; for proofs see 
\cite{Nijenhuis, Richardson, 1967}, 
\cite{Michor, 1987}, and \cite{Lecomte, Michor, Schicke\-tanz} for 
multigraded versions; the extension to the convenient setting does 
not offer any difficulties.
\roster
\item [3] For $K\in\La^{k+1}(V;V)$, $\ph\in\La^p(V)$, and 
       $\Ph\in\La(V;W)$ we have 
       $i_K(\ph\wedge\Ph)=i_K\ph\wedge\Ph+(-1)^{kp}\ph\wedge i_K\Ph$ 
       so $i_K$ is a graded derivation of degree $k$ of the 
       $\La(V)$-module $\La(V;W)$ and any derivation is of that form.
\item The space of graded derivations of the graded $\La(V)$-module 
       $\La(V;W)$ is a graded Lie algebra with bracket the graded 
       commutator $[D_1,D_2]=D_1D_2-(-1)^{d_1d_2}D_2D_1$, see 
       \nmb!{3.1}.
\item For $K\in\La^{k+1}(V)$ and $L\in\La^{\ell+1}(V)$ we have 
	  $[i_K,i_L]=i([K,L]^{\wedge})$ where 
       $[K,L]^{\wedge}=i_KL-(-1)^{k\ell}i_LK$. So by \therosteritem4 
       we get a graded Lie algebra \newline
       $(\La^{*+1}(V;V),[\quad,\quad]^{\wedge})$, called the 
       Nijenhuis-Richardson algebra.
\item If $\mu\in\La^2(V;V)$, i\. e\. $\mu:V\x V\to V$ is bounded skew 
       symmetric bilinear, then $[\mu,\mu]^{\wedge}=2i_\mu\mu=0$ if 
       and only if $(V,\mu)$ is a convenient Lie algebra.
\endroster

\subheading{\nmb.{7.3}. Polyderivations } Let $A$ be a convenient 
algebra and let $L^k(A)\subset \La^{k+1}(A;A)$ be the space of all 
maps $K$ such that for any $a_1,\dots a_k\in A$ the linear map 
$a\mapsto K(a,a_1,\dots ,a_k)$ is a derivation of $A$. Obviously this 
is a closed linear subspace and thus each $L^k(A)$ is a convenient 
vector space. We call $L(A):=\bigoplus _{k\geq 0}L^k(A)$ the space of 
all \idx{\it skew symmetric polyderivations} of $A$.
Obviously $L^k(A)$ is not an $A$ submodule of $\La^{k+1}(A;A)$ in 
general.

\proclaim{\nmb.{7.4}. Theorem } Let $A$ be a convenient algebra. 
Then $(L(A),[\quad ,\quad ]^{\wedge})$ is a graded Lie subalgebra of the 
Nijenhuis-Richardson algebra 
$(\La^{*+1}(A;A),[\quad ,\quad ]^{\wedge})$. 

So $(L(A),[\quad ,\quad ]^{\wedge})$
is a convenient graded Lie algebra called the \idx{\it
Schouten-Ni\-jen\-huis algebra\/} of $A$.
\endproclaim

\demo{Proof}
It suffices to show that for $K_i\in L^{k_i}(A)$ the bracket
$[K_1,K_2]^{\wedge}$ again lies in $L(A)$. This means that we have to show 
that for arbitrary elements $a,b\in A$ we have:
$$i _{ab}[K_1,K_2]^\wedge=(i _a[K_1,K_2]^\wedge)b+a(i _b[K_1,K_2]^\wedge)$$

 From \nmb!{7.2}.(5) we see that for $a\in A = \La^0(A;A)$ and 
$K\in \La^{k+1}(A)$ we have
$$i_ai_K-(-1)^ki_Ki_a=i([a,K]^\wedge)=i(i_aK). \tag 1$$
If furthermore $L\in L^\ell$  we obviously have from the polyderivation 
property of $L$:
$$\gather
i(K\wedge a)L = i_KL\wedge a + K\wedge i_aL, \tag2\\
i(a\wedge K)L = a\wedge i_KL + (-1)^{(k+1)\ell}i_aL\wedge K. \tag3
\endgather$$
Using this we may compute as follows, where we delete $\wedge$ if one 
of the factors is in the algebra $A$:
$$\align
i _{ab}&[K_1,K_2]^\wedge=i _{ab}(i(K_1)K_2)-(-1)^{k_1k_2}i _{ab}
(i(K_2)K_1)=\\
=&i(i _{ab}K_1)K_2+(-1)^{k_1}i(K_1)(i _{ab}K_2)-\\
&-(-1)^{k_1k_2}i(i _{ab}K_2)K_1-(-1)^{(k_1+1)k_2}i(K_2)
(i _{ab}K_1)=\\
=&i(i _aK_1b)K_2)+i(a(i _bK_1)K_2)+\\
&+(-1)^{k_1}i(K_1)(i _aK_2b)+(-1)^{k_1}i(K_1)(ai _bK_2)-\\
&-(-1)^{k_1k_2}i(i _aK_2b)K_1-(-1)^{k_1k_2}i(ai _bK_2)K_1-\\
&-(-1)^{(k_1+1)k_2}i(K_2)(i _aK_1b)-(-1)^{(k_1+1)k_2}
i(K_2)(ai _bK_1)=\\
=&(i(i _aK_1)K_2)b+(i _aK_1)\wedge(i _bK_2)+\\
&+a(i(i _bK_1)K_2)+(-1)^{k_1k_2}(i _aK_2)\wedge(i _bK_1)+\\
&+(-1)^{k_1}(i(K_1)(i _aK_2))b+(-1)^{k_1}a(i(K_1)(i _bK_2))-\\
&-(-1)^{k_1k_2}(i(i _aK_2)K_1)b-(-1)^{k_1k_2}(i _aK_2)\wedge
(i _bK_1)-\\
&-(-1)^{k_1k_2}a(i(i _bK_2)K_1)-(i _aK_1)\wedge(i _bK_2)-\\
&-(-1)^{(k_1+1)k_2}(i(K_2)(i _aK_1))b-(-1)^{(k_1+1)k_2}
a(i(K_2)(i _bK_1))=\\
=&(i _a(i(K_1)K_2))b-(-1)^{k_1k_2}(i _a(i(K_2)K_1))b+\\
&+a(i _b(i(K_1)K_2))-(-1)^{k_1k_2}a(i _b(i(K_2)K_1))=\\
=&(i _a[K_1,K_2]^\wedge)b+a(i _b[K_1,K_2]^\wedge)\qed
\endalign$$
\enddemo

\subheading{\nmb.{7.5}. Definition} Let $A$ be an algebra. A 
2-derivation $\mu\in L^1(A)$ is called a \idx{\it Poisson structure} 
on $A$ if $[\mu,\mu]^{\wedge}=0$. 

\proclaim{\nmb.{7.6}. Theorem} Let $\mu$ be a Poisson structure for 
the algebra $A$. Then $\mu:A\x A\to A$ is a Lie algebra structure.
Furthermore we have  
$$\align
\mu(ab,c)&=a\mu(b,c)+\mu(a,c)b,\\
\mu(a,bc)&=b\mu(a,c)+\mu(a,b)c.
\endalign$$
The mapping $\check\mu: A\to \Der(A), a\mapsto \mu(a,\quad)$ is a 
homomorphism of Lie algebras $(A,\mu)\to (\Der(A),[\quad,\quad])$, 
where the second bracket is the Lie bracket (commutator), see 
\nmb!{4.2}. 
\endproclaim

This is the non-commutative generalization of the Poisson bracket of 
differential geometry.

\demo{Proof}
\nmb!{7.2}.(6) implies that $\mu$ is a Lie algebra structure.
The other assertion is just the property of a polyderivation.
\qed\enddemo


\Refs

\ref
\by Boman, J\.
\paper Differentiability of a function and of its compositions with 
functions of one variable
\jour Math\. Scand\.
\vol 20
\pages 249--268
\yr 1967
\endref

\ref
\by Cartan, H\.; Eilenberg, S\.
\book Homological Algebra
\publ Princeton University Press
\publaddr Princeton
\yr 1956
\endref

\ref
\by Connes, A.
\paper Non-commutative differential geometry
\jour Publ. Math. I.H.E.S.
\vol 62
\pages 257--360
\yr 1985
\endref


\ref
\by Coquereaux, R\.; Kastler, D\.
\paper Remarks on the differential envelopes of associative algebras
\jour Pacific J\. of Math\. 
\vol 137(2)
\yr 1989
\endref

\ref
\by De Wilde, M.; Lecomte, P\. B\. A\. 
\paper Formal deformations of the Poisson Lie algebra of a symplectic manifold and star-products. Existence, equivalence, derivations
\inbook Deformation theory of algebras and structures and applications, M. Hazewinkel, M. Gerstenhaber, Eds
\publ Kluwer Acad. Publ.
\publaddr Dordrecht
\pages 897--960
\yr 1988
\endref

\ref 
\by De Wilde, M.; Lecomte, P\. B\. A\. 
\paper Existence of star products on exact symplectic manifolds
\jour Ann. Inst. Fourier
\vol 35
\pages 117--143
\yr 1985
\endref

\ref 
\by Dubois-Violette, Michel
\paper Derivations et calcul differentiel non-commutatif
\jour C\.R\. Acad\. Sci\. Paris, S\'erie I
\vol 297
\page 403--408
\yr 1988
\endref

\ref \by Fr\"olicher, Alfred; Kriegl, Andreas
\book Linear spaces and differentiation theory 
\bookinfo Pure and Applied Mathematics
\publ J. Wiley
\publaddr Chichester
\yr 1988
\endref

\ref   
\by Fr\"olicher, A.; Nijenhuis, A.   
\paper Theory of vector valued differential forms. Part I  
\jour Indagationes Math   
\vol 18   
\yr 1956   
\pages 338--359   
\endref

\ref
\by Hochschild, G.; Kostant, B.; Rosenberg, A.
\paper Differential forms on regular affine algebras
\jour Trans. Amer. Math. Soc.
\vol 102
\yr 1962
\pages 383--408
\endref

\ref 
\by Jadczyk, A,; Kastler, D.
\paper Graded Lie Cartan pairs I
\jour Rep\. Math\. Phys\.
\vol 25
\pages 1--51
\yr 1987
\endref

\ref 
\by Jadczyk, A,; Kastler, D.
\paper Graded Lie Cartan pairs II: The Fermionic differential calculus
\jour Ann\. Phys\.
\vol 179
\pages 169--200
\yr 1987
\endref

\ref 
\by K\"ahler, E.
\paper Algebra und Differentialrechnung
\inbook Bericht \"uber die Mathematikertagung in Berlin 1953
\endref

\ref   
\by Kainz, G\.; Kriegl, A\.; Michor, P\. W\.   
\paper $C^\infty$-algebras from the functional analytic  viewpoint   
\jour J. pure appl. Algebra   
\vol 46   
\yr 1987   
\pages 89-107   
\endref 

\ref 
\by Karoubi, Max
\paper Connexions, courbures et classes caract\'eristiques en K-theorie alg\'ebriques
\inbook Canadian Math. Soc. Conference Proc. Vol 2
\pages 19--27
\yr 1982
\endref

\ref 
\by Karoubi, Max
\paper Homologie cyclique des groupes et alg\'ebres
\jour C. R. Acad. Sci. Paris
\vol 297
\pages 381--384
\yr 1983
\endref

\ref 
\by Karoubi, Max
\book Homologie cyclique et K-th\'eorie
\bookinfo Asterisque 149
\publ Soci\'et\'e Math\'e\-ma\-ti\-que de France
\yr 1987
\endref

\ref 
\by Kastler, D.; Stora, R.
\paper Lie-Cartan pairs
\jour J. Geometry and Physics
\vol 2
\pages 1--31
\yr 1985
\endref

\ref
\by Keller, H\.H\.
\book Differential calculus in locally convex spaces
\bookinfo Springer Lecture Notes in Mathematics 417
\yr 1974
\endref

\ref 
\by Krasil'shchik, I. S.
\paper Schouten bracket and canonical algebras
\inbook Global Analysis - Studies and Applications III 
\bookinfo Springer Lecture Notes in Mathematics 1334
\pages 79--110
\yr 1988
\endref

\ref
\by Kriegl, Andreas
\paper Die richtigen R\"aume f\"ur Analysis im Unendlich - Dimensionalen
\jour Monatshefte Math.
\vol 94
\yr 1982 
\pages 109--124
\endref

\ref
\by Kriegl, A.
\paper Eine kartesisch abgeschlossene Kategorie glatter Abbildungen
zwischen beliebigen lokal\-konvexen Vektorr\"aumen
\jour Monatshefte f\"ur Math.
\vol 95
\yr 1983
\pages 287--309
\endref

\ref  
\by Kriegl, Andreas; Michor, Peter W.  
\paper A convenient setting for real analytic mappings 
\jour Acta Mathematica 
\vol 165
\pages 105--159
\yr 1990 
\endref

\ref 
\by Kriegl, Andreas; Michor, Peter W.
\book Foundations of Global Analysis  
\bookinfo in preparation  
\endref

\ref
\by Kriegl, Andreas; Nel, Louis D.
\paper A convenient setting for holomorphy
\jour Cahiers Top. G\'eo. Diff.
\vol 26
\yr 1985
\pages 273--309
\endref

\ref 
\by Kunz, Ernst
\book K\"ahler Differentials
\publ Viehweg
\publaddr Braunschweig - Wiesbaden
\yr 1986
\endref

\ref 
\by Lecomte, P. B. A.
\paper Sur la suite exacte canonique associ\'ee a un fibre principal
\jour Bull. Soc. math. France
\vol 113
\pages  259--271
\yr 1985
\endref

\ref 
\by Lecomte, P. B. A.; Melotte, D\.; Roger, C\.
\paper Explicit form and convergence of 1-differential formal deformations of the Poisson Lie algebra
\jour Lett. Math. Physics
\vol 18
\pages 275--285
\yr 1989
\endref

\ref   
\by Lecomte, P\. B\. A\.; Michor, P\. W\.; Schicketanz, H\.   
\paper The multigraded Nijenhuis-Richardson Algebra, its universal property and application
\paperinfo to appear
\jour J. Pure Applied Algebra
\endref

\ref
\by Michor, Peter W\.
\book Manifolds of differentiable mappings
\publ Shiva
\publaddr Orpington
\yr 1980
\endref

\ref   
\by Michor, Peter W.   
\paper Remarks on the Fr\"olicher-Nijenhuis bracket   
\inbook Proceedings of the Conference on Differential Geometry and its Applications, Brno 1986   
\publ D.~Reidel   
\yr 1987   
\pages 197--220   
\endref

\ref   
\by Michor, P. W.   
\paper Knit products of graded Lie algebras and groups    
\paperinfo Proceedings of the Winter School on Geometry and Physics, Srni 1989   
\jour Suppl. Rendiconti Circolo Mat. Palermo, Serie II   
\yr 1990   
\endref

\ref
\by Milnor, John
\paper Remarks on infinite dimensional Lie groups
\inbook Relativity, Groups, and Topology II, Les Houches, 1983,
B\.S\. DeWitt, R\. Stora, Eds\.
\publ Elsevier
\publaddr Amsterdam
\yr 1984
\endref

\ref
\by Narnhofer, H.; Thirring, W\.; Wicklicky, H\.
\paper Transitivity and ergodicity of Quantum systems
\jour J\. Statistical Physics
\vol 52
\yr 1988
\pages 1097--1112
\endref

\ref 
\by Nijenhuis, A.; Richardson, R. 
\paper
Deformation of Lie algebra structures \jour J. Math. Mech. \vol
17 \yr 1967 \pages 89--105 \endref

\ref
\by Pressley, A\.; Segal, G\.
\book Loop groups
\bookinfo Oxford Mathematical Monographs 
\publ Oxford University Press
\yr 1986
\endref

\endRefs
\enddocument